\documentclass[aap,citesort,MSNbibl,dvips]{arximspdf}
\usepackage{mathbh}
\usepackage{dcolumn}
\usepackage{graphicx}

\doi{10.1214/09-AAP659}
\volume{20}
\issue{4}
\pubyear{2010}
\firstpage{1389}
\lastpage{1424}

\makeatletter

\newcolumntype{d}[1]{D{.}{.}{#1}}

\newcommand{\vd}{{d}}

\newtheorem{proposition}{Proposition}

\newproclaim{algorithm}{Algorithm}
\newproclaim{remark}{Remark}
\newproclaim{notation}{Notation}

\newtheorem{lemma}{Lemma}

\makeatother

\begin{document}
\begin{frontmatter}

\title{Simulation of diffusions by means of importance sampling paradigm}
\runtitle{Simulation of diffusions}

\begin{aug}
\author[A]{\fnms{Madalina} \snm{Deaconu}\corref{}\ead[label=e1]{Madalina.Deaconu@inria.fr}} and
\author[A]{\fnms{Antoine} \snm{Lejay}\thanksref{t1}\ead[label=e2]{Antoine.Lejay@iecn.u-nancy.fr}}
\runauthor{M. Deaconu and A. Lejay}
\affiliation{{\'E}quipe Projet TOSCA---Institut Elie Cartan UMR 7502
Nancy-Universit{\'e},~CNRS,~INRIA}
\address[A]{IECN\\
Campus Scientifique\\
BP 70239\\
54506 Vand\oe uvre-l{\`e}s-Nancy cedex\\
France\\
\printead{e1}\\
\phantom{E-mail: }\printead*{e2}} 
\end{aug}

\thankstext{t1}{Supported in part by the GdR MOMAS
(funded by ANDRA, BRGM, CEA, CNRS, EDF and IRSN).}

\received{\smonth{1} \syear{2007}}
\revised{\smonth{8} \syear{2009}}

%
\begin{abstract}
The aim of this paper is to introduce a new Monte Carlo method based on
importance sampling techniques for the simulation of stochastic
differential equations. The main idea  is to combine random walk on
squares or rectangles methods with importance sampling techniques.

The first interest of this approach is that the weights can be easily
computed from the density of the one-dimensional Brownian motion.
Compared to the Euler scheme this method allows one to obtain a more
accurate approximation of diffusions when one has to consider complex
boundary conditions. The method provides also an interesting
alternative to performing variance reduction techniques and simulating
rare events.
\end{abstract}

%

\end{frontmatter}

\section{Introduction}\label{sec1}

Monte Carlo methods are sometimes the unique alternative used to
solve numerically partially differential equations
(PDE) involving an operator
of the form
\[
L=\frac{1}{2}\sum_{i,j=1}^d a_{i,j}(\cdot)\frac{\partial
^2}{\partial
x_i\,\partial x_j}
+\sum_{i=1}^d b_i(\cdot)\,\frac{\partial}{\partial x_i}.
\]
The operator $L$ is the infinitesimal generator associated with the solution
of the stochastic differential equation (SDE)
%
%
\begin{equation}
\label{SDE}
X_t=X_0+\int_0^t \sigma(X_s)\,\vd B_s
+\int_0^t b(X_s)\,\vd s\qquad\mbox{with }\sigma{\sigma}^{*}=a.
\end{equation}
It is well known that, for $T>0$ fixed, the solution on the cylinder
$[0,T]\times D$, of the parabolic PDE,
\[
\cases{
\dfrac{\partial u(t,x)}{\partial t}+Lu(t,x)=0,\vspace*{2pt}\cr
u(T,x)=g(x), &\quad for $x \in D$,\cr
u(t,x)=\phi(t,x), &\quad for $(t,x)\in[0,T]\times\partial D$,}
\]
can be written as
\[
u(t,x)=\mathbb{E}_{t,x}[g(X_T);T\leq\tau]
+\mathbb{E}_{t,x}[\phi(\tau,X_\tau);\tau<T],
\]
where $\tau$ stands for the first exit time of $X$ from
the domain $D$. $\mathbb{E}_{t,x}$ means that the process $X$ is
starting from $x$ at time $t$. Thus an approximation of $u(t,x)$ can be obtained
by averaging $g(X_T)\mathbh{1}_{T\leq\tau}$ and $\phi(\tau,X_\tau
)\mathbh{1}_{\tau<T}$
over a large number of realizations of paths
of $X$. Elliptic
PDE may be considered as well.

A large spectra of methods has been already proposed in order
to simulate $X$
(see, e.g., the books of Kloeden and Platen \cite{kloeden92a}
and of Milstein and Tretyakov \cite{milstein04a}).
Most of these methods are extensions of the
Euler scheme which provides a very efficient way
to simulate (\ref{SDE}) in the whole space. This method
becomes harder to set up in a bounded domain,
either with an absorbing or a reflecting boundary
condition. Nevertheless some refinements have been proposed
(see, e.g.,
\cite{bossy04a,gobet00a,gobet01a,jansons03a,pettersson95a,slominski01a}).
To improve the quality of the simulation or to speed it up, variance
reduction techniques can be considered (see, e.g.,
\cite
{arouna04a,arouna04b,bardou05a,kebaier05a,heath02a,kohatsu02a,newton94a,zou04a}).
This list is not intended to be exhaustive.

In the simplest situation, for $a=\mathrm{Id}$ and $b=0$, the
underlying diffusion process is the Brownian motion. Muller
proposed in 1956 a very simple scheme to solve a Dirichlet
boundary value problem. This method is called
the \textit{random walk on spheres method}
\cite{muller56a}. The idea is to simulate successively,
for the Brownian motion, the
first exit position from the largest sphere included in the domain and
centered in the starting point. This exit position becomes the new
starting point, and the procedure is iterated until the exit
point is close enough to the
boundary. Nevertheless, simulating the exit time from a sphere is
numerically costly. In~\cite{milstein93a},
Milstein and Rybkina proposed to use this scheme for
solving (\ref{SDE}) by freezing locally the value of the
coefficients. In a first approach, spheres (that become
ellipsoids) were used. Later on
\cite{milstein95a} (see also \cite{milstein04a}),
Milstein and Tretyakov used
time--space parallelepipeds
with a cubic space basis. For this last
approach, it is easier to keep track of the time but the
involved random variables are costly to simulate. In order to overcome
these difficulties,
one may think to use tabulated values. This is memory
consuming as the random variables to simulate depend on
one or two parameters. The method of random walk
on squares was also independently developed in
the Ph.D. thesis of Faure \cite{faure92a}.
For the Brownian motion, this
method is still a good alternative to the random walk on
spheres (see \cite{lejaym} for an application in
geophysics).

In \cite{deaconulejay05a}, we have proposed a scheme for
simulating the exact exit time and position from a rectangle
for the Brownian motion starting from any point inside this rectangle.
Compared to the random walk on spheres method,
this method has the following
advantages:
\begin{itemize}
\item It can be used
whatever the dimension and, as for the random walk on squares,
a constant drift term may be added.
\item The rectangles can be chosen prior to any simulation, and not
dynamically. There is no need to consider smaller and
smaller spheres or squares when the particle is near the boundary.
\item The method can be also adapted and used for the simulation
of diffusion processes killed on some part of the boundary.
\end{itemize}

The method we propose here is based on the idea to
simulate the first exit time and position from a
parallelepiped by using an importance sampling technique
(see, e.g., \cite{fishman96a,glasserman04a}). The
exit time and position from a parallelepiped for a Brownian
motion with locally frozen coefficients is chosen
arbitrarily, and a weight is computed at each simulation. By
repeating this procedure, we
get the density on the boundary or at a given time of
the particles, by weighting the simulated paths.
As we will see, the
weights are rather easily deduced from the density of the
one-dimensional Brownian motion killed when it exits from
$[-1,1]$. All involved expressions are numerically easy
to implement.

This new algorithm is slower than the Euler scheme for smooth
coefficients, but it is faster than
the random walk on squares \cite{lejaym,milstein04a}
and the random walk on rectangles \cite{deaconulejay05a}.
It can be used to simulate the Brownian motion as well as
solutions of stochastic differential equations for specific
complex situations as:
(a) complex geometries (the boundary conditions are correctly taken into
account);
(b) fast estimation of the exit time of a domain for the Brownian
motion (only few rectangles are needed);
(c) variance reduction;
(d) simulation of rare events.

This algorithm could be relevant for many domains:
finance, physics, biology, geophysics, etc. It may also be
used locally (e.g., it can be mixed with the Euler
scheme and used when the particle is close to the
boundary) or combined with other algorithms, such as
population Monte Carlo methods (see Section \ref{sec-population}).

We conclude this article with numerical simulations
illustrating various examples. It has to be
noted that choosing ``good'' distributions for the exit time
and position from a rectangle is not an easy task in order
to reduce the variance. We then plan to study in the future
how to construct algorithms that minimize the variance,
as in \cite{arouna04a,bardou05a}. We have to consider for this
a high-dimension optimization problem.

\subsection*{Outline} In Section \ref{sec-para},
we present the importance sampling technique applied
to the exit time and position for a (drifted) Brownian
motion from a rectangle. In Section \ref{sec-density},
we recall briefly some results about the density
of the one-dimensional Brownian motion with different
boundary conditions. The explicit expressions
are given in the \hyperref[app-1]{Appendix}. In Section \ref{sec-general},
we present our algorithm and compute its weak error.
Four test cases are presented in Section \ref{sec-numerical}. We
compare also our algorithm with other methods in this last
section.


\section{Algorithm for the exit time and position from a right
time--space parallelepiped by using an importance sampling method}
\label{sec-para}

The aim of this part is to give a clear presentation of our method.
In order to avoid ambiguous notation we consider in this section
the situation of a two-dimensional space domain. The results can
be easily generalized to higher space dimension.

We are looking for an accurate approximation of the exit time
and position from a right time--space parallelepiped
which is a geometric figure in the three-dimensional space.

For $L_1,L_2>0$ given let $R$ be
the rectangle $[-L_1,L_1]\times[-L_2,L_2]$. The rectangle $R$
is the space basis of the right time--space parallelepiped
$R_T=[0,T]\times R$
for a fixed $T>0$. We can also consider $R_\infty=\mathbb{R}_+\times R$,
and set in this case $T=+\infty$.

\begin{figure}[b]

\includegraphics{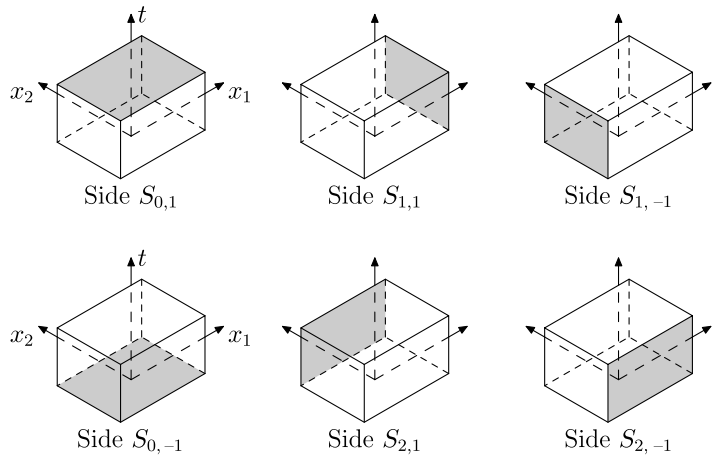}

\caption{Convention for the sides of $R_T=[0,T]\times R$.}
\label{fig-1}
\end{figure}

For $T<+\infty$, the right time--space parallelepiped $R_T$
has six sides which are denoted by
\begin{eqnarray*}
S_{0,1}&=&\{T\}\times R,\\
S_{0,-1}&=&\{0\}\times R,\\
S_{1,\eta}&=&[0,T]\times[-L_1,L_1]\times\{\eta L_2\} \qquad\mbox{for
}\eta\in\{-1,1\},\\
S_{2,\eta}&=&[0,T]\times\{\eta L_1\}\times[-L_2,L_2] \qquad\mbox{for
}\eta\in\{-1,1\}.
\end{eqnarray*}
In other words, each side of $R_T$ is labeled by
a couple $(i,\eta)\in\{0,1,2\}\times\{-1,1\}$.
For $i\in\{1,2\}$ the side $S_{i,\eta}$ is
perpendicular to the unit vector in the
$i$th direction.
For $i=0$, the side $S_{0,-1}$ corresponds to
the rectangular initial basis while the side $S_{0,1}$ corresponds
to the top of the time--space parallelepiped $R_T$ for $T<+\infty$
(see Figure \ref{fig-1}).

From now on, we shall identify each side
with the corresponding $(i,\eta)$-indices.

%
%

We consider a time-homogeneous diffusion process $(X_t)_{t\geq0}$
living in $R$. On each side of $R$, the process $X$ may be reflected or
absorbed. Moreover, if $T<+\infty$, the process
is stopped at time $T$. We can thus identify the sides of $R$
with the sides $S_{i,\eta}$ of $R_T$ for $i\in\{1,2\}$ and $\eta\in
\{-1,1\}$.
We denote by $\mathfrak{R}$ the subset of $\{1,2\}\times\{-1,1\}$
that contains the indices of the sides
on which a Neumann boundary condition holds (possibly,
$\mathfrak{R}=\varnothing$). On this set
the diffusion is reflected. Let us denote by $\mathfrak{D}$ the subset
of $\{1,2\}\times\{-1,1\}$ that contains the indices
of the sides on which a Dirichlet boundary condition
holds. On this set the diffusion is killed. Finally
let us set $\mathfrak{A}=\mathfrak{D}$ if $T=+\infty$ and $\mathfrak
{A}=\{(0,1)\}\cup\mathfrak{D}$
if $T<+\infty$. With this notation the time--space process
$t\mapsto(t,X_t)$ is absorbed
when hitting one of the sides $S_{i,\eta}$
with $(i,\eta)\in\mathfrak{A}$.

Let $B=(B^1,B^2)$ be a two-dimensional Brownian
motion and $\mu=(\mu_1,\mu_2)$ a vector of $\mathbb{R}^2$.
For $i\in\{1,2\}$, we set
\[
\gamma_{i,\eta}=\cases{
1, &\quad if $(i,\eta)\in\mathfrak{R}\mbox{ (reflection)}$,\cr
0, &\quad if $(i,\eta)\in\mathfrak{A}\mbox{ (absorption)}$.}
\]

We consider the two-dimensional diffusion process
$(X,\mathbb{P}_x)_{x\in R}$ whose coordinates are, for $x=(x_1,x_2)\in R$,
%
%
\begin{equation}
\label{eq-2}
X_t^i=x_i+B^i_t+\mu_i
t+\gamma_{i,1}\ell_t^{L_i}(X^i)-\gamma_{i,-1}\ell
_t^{-L_i}(X^i),\qquad
\mathbb{P}_x\mbox{-a.s.},
\end{equation}
where $\ell_t^{\pm L_i}(X^i)$ stands for the symmetric local time
of $X^i$ at $\pm L_i$, respectively.

We define $\tau_0=T$, $\tau_i=\inf\{t>0 | |X_t^i|>L_i\}$
for $i\in\{1,2\}$ and
\[
\tau= \min_{i\in\{0,1,2\}}\tau_i.
\]

In addition, we set $J=\arg\min_{i\in\{0,1,2\}} \tau_i$. With
this notation, unless $J\not=0$, the $J$th component of $X$ is
the first to exit from the domain.
For $J\in\{1,2\}$, let us define $\varepsilon=X_{\tau_J}^J/L_J\in\{
-1,1\}$.
For $J=0$ we set $\varepsilon=1$. In this case $X$ has
not reached the sides of $\mathfrak{D}$ before time $T$.

The couple $(J,\varepsilon)$ labels the side in $\mathfrak{A}$
of the parallelepiped $R_T=[0,T]\times R$ that the diffusion $X$ hits first.
Note that with our convention, the sides on which the process is reflected
cannot be reached so that $\tau_i=+\infty$
if $X^i$ is reflected both at $-L_i$ and $L_i$.

We are interested
in computing $\mathbb{E}_x[f(\tau,X_\tau)]$
by a Monte Carlo method for a bounded, measurable function $f$ where
$\tau$ is defined as above.

Instead of simulating $(\tau,X_\tau)$,
we will simulate some random variables according
to the following procedure. The aim is
to simulate $(J,\varepsilon,\tau,X_\tau)$ by
using an importance sampling technique.
In order to do this we choose a probability $\widehat{\mathbb{P}}_x$
which is absolutely continuous with respect to $\mathbb{P}_x$,
and we draw a realization of $(J,\varepsilon,\tau,X_\tau)$.
Let us set
\[
\alpha_{i,\eta}=\widehat{\mathbb{P}}_x[(J,\varepsilon)=(i,\eta)]
\]
for $(i,\eta)\in\mathfrak{A}$. For $(i, \eta)\in\mathfrak{A}$ let
$k_{i,\eta}$
denote the density under $\widehat{\mathbb{P}}_x$ of $(\tau,X_\tau)$
given $\{(\tau,X_\tau)\in S_{i,\eta}\}$.

In order to simplify notation let us consider an underlying
probability space $(\Omega, \mathcal{{F}}, \mathbb{P}_x)$ rich enough.
Let $Z$ be a random variable on this space, with distribution~$\mathbb{P}_x$.
Let $A$ be a measurable event on this space.
We suppose that, conditionally on $A$, $Z$ has a density $p(\cdot|A)$
with respect to the Lebesgue measure. Let us introduce the following convention:
\[
\mathbb{P}_x[Z=z; A]=p(z|A)\mathbb{P}_x[A].
\]
That is, for $B$ a measurable event of $(\Omega, \mathcal{{F}},
\mathbb{P}_x)$,
\[
\mathbb{P}_x[\{Z\in B\}\cap A]
=\int_B p(z|A)\mathbb{P}_x[A]\,\vd z=\int_B \mathbb{P}_x[Z=z; A]\,
\vd z.
\]
%
Consider now the following notation:
let $(i,\eta)\in\mathfrak{A}$. For $i\in\{1,2\}$ set $j=3-i$.
Then for any $\theta>0$ and $z\in S_{i,\eta}$, we define
%
%
\begin{equation}
\label{eq-weight-1}
M_{i,\eta}(\theta,z)=\frac
{\mathbb{P}_{x}[\tau_i=\theta; X^i_{\tau_i}=\eta L_i]
\mathbb{P}_{x}[X^j_\theta=z_j; \tau_j>\theta]}
{\alpha_{i,\eta}k_{i,\eta}(\theta,z)},
\end{equation}
where $k_{i,\eta}$ is the $\{X_\tau\in S_{i,\eta}\}$-conditional
density under $\widehat{\mathbb{P}}_x$ of $(\tau,X_\tau)$.

If $T<+\infty$, we define
%
%
\begin{equation}
\label{eq-weight-2}
M_{0,1}(T,z)
=\frac{1}{\alpha_{0,1}k_{0,1}(T,z)}
\prod_{j\in\{1,2\}}\mathbb{P}_{x}[X^j_T=z_j; \tau_j>T],
\end{equation}
where $k_{i,\eta}$ is the $\{X_\tau\in S_{i,\eta}\}$-conditional
density under $\widehat{\mathbb{P}}_x$ of $(\tau,X_\tau)$.

We call $M_{i,\eta}$ \textit{weight}.
\begin{proposition}\label{prop1}
The \textup{weights} $M_{i,\eta}$ defined in $(\ref{eq-weight-1})$ and
$(\ref{eq-weight-2})$ satisfy
\[
\mathbb{E}_x[f(\tau,X_\tau)]
=\widehat{\mathbb{E}}_x[M_{J,\varepsilon}(\tau,X_{\tau})f(\tau
,X_{\tau})]
\]
for any measurable and bounded function $f$ on $\partial R_T$.
\end{proposition}

Before proving this proposition let us introduce the algorithm.

The algorithm is described as follows:
\begin{algorithm}
\label{algo-1}
Let $x$ be fixed in $R$.
\begin{enumerate}[(1)]
\item[(1)] Draw a realization $(\overline{J},\overline{\varepsilon})$
of $(J,\varepsilon)\in\mathfrak{A}$ under $\widehat{\mathbb{P}}_x$.
\item[(2)] Draw a realization of the exit time and exit position
$(\overline{\tau},\overline{X}_{\overline{\tau}})$
according to the density $k_{\overline{J},\overline{\varepsilon}}$ on
$S_{\overline{J},\overline{\varepsilon}}$.
\item[(3)] Compute
the value of $M_{\overline{J},\overline{\varepsilon}}(\overline
{\tau},\overline{X}_{\overline{\tau}})$ by
\[
\widehat{\mathbb{E}}_x[M_{J,\varepsilon}(\tau,X_\tau)f(\tau
,X_\tau)]=\mathbb{E}_x[f(\tau,X_\tau)].
\]
We call $M_{\overline{J},\overline{\varepsilon}}(\overline{\tau
},\overline{X}_{\overline{\tau}})$, \textit{weight}.
\end{enumerate}
\end{algorithm}

If
$\{(\overline{J}^{(i)},\overline{\varepsilon}^{(i)},\overline{\tau
}^{(i)},\overline{X}_{\overline{\tau}}^{(i)},\overline{w}^{(i)})\}
_{i=1,\ldots,N}$
are $N$ independent realizations of the random variables
$(J,\varepsilon,\tau,X_\tau,M_{J,\varepsilon}(\tau,X_\tau))$
constructed as above, by
the law of large numbers we have
\[
\mathbb{E}_x[f(\tau,X_\tau)]
=\lim_{N\to\infty}\frac{1}{N}\sum_{i=1}^N
\overline{w}^{(i)}f\bigl(\overline{\tau}^{(i)},\overline
{X}_{\overline
{\tau}^{(i)}}^{(i)}\bigr).
\]

The main feature of our approach is that the weights
$M_{J,\varepsilon}(\tau,X_\tau)$
can be easily evaluated.
\begin{remark}
In order to evaluate $M_{i,\eta}$ with (\ref{eq-weight-1})
and (\ref{eq-weight-2}),
there is no need to know $\mathbb{P}_x[(J,\varepsilon)=(i,\eta)]$. It
is important to notice that $M_{i,\eta}$ depends only on the
one-dimensional distributions
of the drifted Brownian motion.
\end{remark}
\begin{pf*}{Proof of the Proposition \protect\ref{prop1}}
We want to prove that
\[
\mathbb{E}_x[f(\tau,X_\tau)]
=\widehat{\mathbb{E}}_x[M_{J,\varepsilon}(\tau,X_{\tau})f(\tau
,X_{\tau})]
\]
for any measurable and bounded function $f$ on $\partial R_T$.

We remark first that if $p_{i,\eta}=\mathbb{P}_x[(J,\varepsilon
)=(i,\eta)]$
for $(i,\eta)$ in $\mathfrak{A}$,
then
%
%
\begin{equation}
\label{eq-23}
\mathbb{E}_x[f(\tau,X_\tau)]
=
\sum_{(i,\eta)\in\mathfrak{A}}
\frac{p_{i,\eta}}{\alpha_{i,\eta}}
\widehat{\mathbb{E}}_x[M_{i,\eta}(\tau,X_\tau)f(\tau,X_\tau) |
(J,\varepsilon)=(i,\eta)].
\end{equation}
Furthermore, for $(i,\eta)\in\mathfrak{D}$,
if $i=2$ set $j=1$ and $z=(z_1,\eta L_2)$ else, if $i=1$ set $j=2$
and $z=(\eta L_1,z_2)$.
\begin{eqnarray*}
&&\mathbb{E}_x[f(\tau,X_\tau) | (J,\varepsilon)=(i,\eta)]\\
&&\qquad=\int_{[0,T]\times[-L_j,L_j]}
f(\theta,z)\mathbb{P}_x[(\tau_i,X^j_{\tau_i})=(\theta,z_j) |
(J,\varepsilon)=(i,\eta)]
\,\vd\theta\,\vd z_j,
\end{eqnarray*}
where $\mathbb{P}_x[(\tau_i,X^j_{\tau_i})=(\theta,z_j) |
(J,\varepsilon)=(i,\eta)]$
is the
$\{(J,\varepsilon)=(i,\eta)\}$-conditional density of $(\tau
_i,X^j_{\tau_i})$
with respect to $\vd t\,\vd z_j$. Hence
\[
\mathbb{E}_x[f(\tau,X_\tau) | (J,\varepsilon)=(i,\eta)]
=\widehat{\mathbb{E}}_x[f(\tau,X_{\tau})M'_{i,\eta}(\tau,X_{\tau
}) | (J,\varepsilon)=(i,\eta) ],
\]
where
\[
M_{i,\eta}'(\theta,z)=
\frac{\mathbb{P}_x[(\tau_i,X^j_{\tau_i})=(\theta,z_j) |
(J,\varepsilon)=(i,\eta)]}{
k_{i,\eta}(\tau,X_{\tau})}.
\]
Let us note that $M_{i,\eta}(\theta,z)=M'_{i,\eta}(\theta
,z)p_{i,\eta}/\alpha_{i,\eta}$. With (\ref{eq-23}), we can deduce that
\[
\mathbb{E}_x[f(\tau,X_\tau)]
=\widehat{\mathbb{E}}_x[f(\tau,X_{\tau})M_{J,\varepsilon}(\tau
,X_{\tau})].
\]

Indeed, it suffices to remark that for $(i,\eta)\in\mathfrak{D}$,
\begin{eqnarray*}
M_{i,\eta}(\theta,z)&=&\frac{1}{\alpha_{i,\eta}k_{i,\eta}(\theta,z)}
\mathbb{P}_x[(\tau_i,X^j_{\tau_i})=(\theta,z_j); (J,\varepsilon
)=(i,\eta)]\\
&=&\frac{1}{\alpha_{i,\eta}k_{i,\eta}(\theta,z)}
\mathbb{P}_x[(\tau_i,X^j_{\tau_i})=(\theta,z_j); X^i_{\tau_i}=\eta
L_i, \tau^j>\theta].
\end{eqnarray*}
The independence of the coordinates of $X$ leads to the desired equality.
If $T<+\infty$, similar computations imply that for $z\in R$,
\[
M_{0,1}(T,z)=
\frac{1}{\alpha_{0,1}k_{0,1}(T,z)}
\mathbb{P}_{x}\Bigl[X_T=z; \min_{i\in\{1,2\}}\tau_i>T \Bigr]
\]
and the conclusion also holds.
\end{pf*}

Let us evaluate these probabilities.

For $i\in\{1,2\}$,
let $p^i(t,x_1,x_2)$ be the solution of
%
%
\begin{equation}
\label{eq-fs}
\cases{
\displaystyle\frac{\partial p^i(t,x_1,x_2)}{\partial t}
=\frac{1}{2}\,\frac{\partial^2 p^i}{\partial x_2^2}(t,x_1,x_2)
+\mu_i\,\frac{\partial p^i}{\partial x_2}(t,x_1,x_2),\cr
\hspace*{134.2pt}\mbox{for }\displaystyle(t,x_1,x_2)\in\mathbb
{R}_+\times(-L_i,L_i)^2,\cr
p^i(t,x_1,x_2)\mathop{\longrightarrow}\limits_{t\searrow0}\delta
_{x_1}(x_2), \qquad\mbox{for $(x_1,x_2)\in(-L_i,L_i)^2$},}
\end{equation}
with the following boundary conditions (b.c.):
\begin{eqnarray*}
p^i(t,x_1,-L_i)&=&0\mbox{ (Dirichlet b.c.)}\qquad\mbox{if }(i,-1)\in
\mathfrak
{A},\\
\frac{\partial p^i}{\partial x_2}(t,x_1,-L_i)&=&0\mbox
{ (Neumann b.c.)}\qquad
\mbox{if }(i,-1)\in\mathfrak{R},\\
p^i(t,x_1,L_i)&=&0\mbox{ (Dirichlet b.c.)}\qquad\mbox{if }(i,1)\in\mathfrak
{A},\\
\frac{\partial p^i}{\partial x_2}(t,x_1,L_i)&=&0\mbox{ (Neumann b.c.)}
\qquad\mbox{if }(i,1)\in\mathfrak{R}.
\end{eqnarray*}
Thus, $p^i$ denotes the density of the drifted Brownian
motion $X^i$ with possibly some reflection at the endpoints
of $(-L_i,L_i)$, and killed when it exits from this
interval by an endpoint where no reflection holds.
For $f$ a bounded measurable function from $[-L_i,L_i]$
to $\mathbb{R}$, we have
\[
\mathbb{E}_{x_1}[f(X^i_t);t<\tau_i]=\int_{-L_i}^{L_i}
p^i(t,x_1,x_2)f(x_2)\,\vd x_2
\]
for $x_1\in[-L_i,L_i]$ where $\mathbb{P}_{x_1}$ is the distribution
of $X^i$
with $X^i_0=x_1\in[-L_i,L_i]$. Let us note that the distribution
of the marginal $X^i$ of $X$ under $\mathbb{P}_{(x_1,x_2)}$ depends
only on $x_i$.

We introduce the \textit{scale function} $\Phi^{i,+}$
of $X^i$ defined by
\[
\mbox{for $x_2\in[-L_i,L_i]$,}\qquad
\Phi^{i,+}(x_2)=
\cases{\displaystyle\frac{e^{2\mu_i L_i}-e^{-2\mu_i
x_2}}{e^{2\mu_iL_i}-e^{-2\mu_iL_i}},
&\quad if $\mu_i\not=0$,\cr
\displaystyle\frac{x_2+L_i}{2L_i}, &\quad if $\mu_i=0$.}
\]
The function $\Phi^{i,+}(x_2)$ has been normalized
such that $\Phi^{i,+}(L_i)=1$.
Let us note that
$\Phi^{i,+}(x_i)=\mathbb{P}_{x_i}[X^i_{\tau_i}=L_i]$
if Dirichlet boundary conditions hold at both endpoints $-L_i$ and $L_i$.
We also set $\Phi^{i,-}(x_2)=1-\Phi^{i,+}(x_2)$.

If Dirichlet boundary conditions hold both at
$-L_i$ and $L_i$, then
we set for $t>0$ and $(x_1,x_2)\in[-L_i,L_i]^2$,
\[
p^{i,\pm}(t,x_1,x_2)=p^i(t,x_1,x_2)\frac{\Phi^{i,\pm}(x_2)}{\Phi
^{i,\pm}(x_1)}.
\]
Via a Doob transform, for a bounded and
measurable function $f$,
\[
\mathbb{E}_{x_1}[f(X_t^i);t<\tau_i | X^i_{\tau_i}=\pm L_i]
=\int_{-L_i}^{L_i} p^{i,\pm}(t,x_1,x_2)f(x_2)\,\vd x_2.
\]

Let us set for $x_1\in(-L_i,L_i)$,
%
%
\begin{equation}
\label{eq-det-1}
q^i(t,x_1)=-\int_{-L_i}^{L_i} \frac{\partial p^i}{\partial
t}(t,x_1,x_2)f(x_2)\,\vd x_2
\end{equation}
and
%
\begin{equation}
\label{eq-det-2}
q^{i,\pm}(t,x_1)=-\int_{-L_i}^{L_i} \frac{\partial p^{i,\pm
}}{\partial t}(t,x_1,x_2)f(x_2)\,\vd x_2.
\end{equation}
We can easily deduce that
\[
\mathbb{P}_{x_1}[\tau_i\leq t]=
\int_0^t q^i(s,x_1)\,\vd s
\quad\mbox{and}\quad
\mathbb{P}_{x_1}[\tau_i\leq t | X^i_{\tau_i}=\pm L_i]=
\int_0^t q^{i,\pm}(s,x_1)\,\vd s.
\]
In other words, $q^i(t,x_1)$ [respectively, $q^{i,\pm}(t,x_1)$]
is the density of the first exit time from $[-L_i,L_i]$
for $X^i$ (respectively, the first exit time from $[-L_i,L_i]$
for $X^i$ given $\{X^i_{\tau_i}=\pm L_i\}$).

Thanks to these expressions, $M_{0,1}(T,z)$ and
$M_{i,\eta}(\theta,z)$ are easily computed since
\begin{eqnarray*}
\mathbb{P}_{x_i}[X^i_\theta=z_i; \tau_i>T]&=&p^i(\theta,x_i,z_i),\\
\mathbb{P}_{x_i}[\tau_i=\theta; X^i_\theta=\pm L_i]
&=&q^{i,\pm}(\theta,x_i)\Phi^{i,\pm}(x_i)\qquad\mbox{if
}(i,-1)\in
\mathfrak{A}
\mbox{ and }(i,1)\in\mathfrak{A},\\
\mathbb{P}_{x_i}[\tau_i=\theta; X^i_\theta=L_i]
&=&q^{i}(\theta,x_i)\qquad\mbox{if }(i,-1)\in\mathfrak{R}\mbox{ and
}(i,1)\in\mathfrak{A},\\
\mathbb{P}_{x_i}[\tau_i=\theta; X^i_\theta=-L_i]
&=&q^{i}(\theta,x_i)\qquad\mbox{if }(i,-1)\in\mathfrak{A}\mbox{ and
}(i,1)\in\mathfrak{R}.
\end{eqnarray*}

\section{Analytical expressions for the densities}
\label{sec-density}

In order to compute $p^i(t,x_1,x_2)$, together with $q^i(t,x_1)$
and $q^{i,\pm}(t,x_1)$ by (\ref{eq-det-1}) and (\ref{eq-det-2}), one has
to solve equation (\ref{eq-fs}). By using a scaling principle, we may
assume that $L_i=1$, as
\[
p^i(t,x_1,x_2)=\frac{1}{L_i}p \biggl(\frac{t}{L_i^2},
\frac{x_1}{L_i},\frac{x_2}{L_i};L_i\mu\biggr),
\]
where $p(t,x_1,x_2;\delta)$ is solution to (\ref{eq-fs})
with $L_i=1$ and a convective term $\mu_i$ equal to $\delta$.

There are basically two ways to obtain $p(t,x_1,x_2;\delta)$.
The first one is based on the spectral
expansion of $\frac{1}{2}\triangle+\delta\nabla$
since this operator may be reduced to
a self-adjoint one with respect
to the scalar product induced by the measure
$\exp(-2\delta x_1)$. The second one is the
method of images when $\delta=0$.

If $\delta\not=0$, the case of a Dirichlet boundary
condition at both endpoints may be treated by using
a simple transform that reduces the problem
to $\delta=0$.

For the case of Neumann boundary condition
at both endpoints, one can invert term by term
the Laplace transform of a series for
the Green function.

In the case of a mixed boundary
condition, the previous method gives rise to a
series that cannot be used in practice, so
only the spectral expansion should be used. In addition,
the first eigenvalues have to be computed numerically.

As the formula are standard in most of the cases,
we give the relevant expressions in
the \hyperref[app-1]{Appendix}.


\section{General domain}
\label{sec-general}

As stated before, we aim to solve by a Monte Carlo method
a parabolic or an elliptic PDE. The idea is
to represent the domain as the union of time--space
parallelepipeds and to simulate the successive
exit times and positions from these parallelepipeds.
Attention has to be paid while doing this decomposition
in order to control the error at each simulation step.

\subsection{From parallelepipeds to right parallelepipeds}
\label{sub-reduc-para}

Consider herein the notation of Section \ref{sec-para}.
Let us study first the parabolic PDE with constant
coefficients
$\lambda$, $c$ and $\mu=(\mu_i)_{i=1,\ldots,d}$
on the rectangle $R_T$,
%
%
\begin{equation}
\label{eq-pde-3}\qquad
\cases{
\displaystyle\frac{\partial v(t,x)}{\partial t}
+\frac{1}{2}\sum_{i=1}^d\frac{\partial^2 v(t,x)}{\partial
x_i^2}\cr
\qquad{}+\displaystyle\sum_{i=1}^d\mu_i\,\frac{\partial
v(t,x)}{\partial
x_i} + cv(t,x)=\lambda, &\quad on
$R_T$,\cr
\dfrac{\partial v(t,x)}{\partial x_i}=0, &\quad for
$x\in S_{i,\eta}\mbox{ if }(i,\eta)\in\mathfrak{R}$,\cr
v(t,x)=\phi(t,x), &\quad
for $x\in S_{i,\eta}\mbox{ if }(i,\eta)\in\mathfrak{A}$,\cr
v(T,x)=g(x), &\quad if $T<+\infty$.}
\end{equation}
We assume that a classical solution to this problem exists,
which is, for example, the case if $\phi$ and $g$ are continuous and bounded.
Let $X$ be the diffusion process whose components
are given by (\ref{eq-2}).
Then it follows from the It{\^o} formula applied to $X$
that, for $t\in[0,T]$,
\begin{eqnarray*}
v(t,x)&=&\mathbb{E}_x\bigl[e^{c(\tau-t)}\phi(\tau-t,X_{\tau
-t});\tau
<T-t\bigr]\\
&&{} +\mathbb{E}_x\bigl[e^{c(T-t)}g(X_{T-t});\tau=T-t\bigr]
+\mathbb{E}_x\biggl[\lambda\int_0^{\tau-t} e^{c(\tau-t-s)}\,\vd s
\biggr],
\end{eqnarray*}
where $\tau$ is as above the first exit time from $R_T$.

Let us remark that if $\sigma$ is an invertible $d\times d$-matrix,
then the function $u(t,x)=v(t,\sigma^{-1} x)$
is solution to
%
%
\begin{equation}
\label{eq-pde-4}
\cases{
\displaystyle\frac{\partial u(t,x)}{\partial t}
+\frac{1}{2}\sum_{i,j=1}^d[\sigma{\sigma}^{*}]_{i,j}\,\frac
{\partial^2 u(t,x)}{\partial x_i\,\partial
x_j}\cr
\qquad{} + \displaystyle\sum_{i=1}^d[\mu{\sigma}^{*}]_i\,\frac
{\partial
u(t,x)}{\partial x_i}
+cu(t,x)=\lambda, \cr
\hspace*{119.08pt}\mbox{on $[0,T]\times\sigma R$},\cr
\sigma_{j,i}\,\dfrac{\partial u(t,x)}{\partial x_j}=0, \qquad\hspace
*{17.7pt}
\mbox{for $x\in\sigma S_{i,\eta}$ if
$(i,\eta)\in\mathfrak{R}$},\cr
u(t,x)=\phi(t,\sigma^{-1}x), \qquad\mbox{for $x\in\sigma S_{i,\eta}$
if $(i,\eta)\in\mathfrak{A}$},\cr
u(T,x)=g(\sigma^{-1}x), \qquad\hspace*{7.09pt} \mbox{if $T<+\infty
$.}}\hspace*{-20pt}
\end{equation}
If $\mathbf{n}_{i}$ is the unit vector orthogonal to the side $\sigma
S_{i,\eta}$,
then $\mathbf{n}_i=({\sigma}^{*})^{-1}\mathbf{e}_i$, where $\mathbf
{e}_i$ is the
unit vector
in the $i$th direction. It follows that $\sigma{\sigma}^{*}\mathbf{n}_i
=\sigma\mathbf{e}_i$ and thus
\[
\mbox{for }x\in\sigma S_{i,\pm1}\qquad
[\sigma{\sigma}^{*}]\mathbf{n}_i\cdot\nabla u(t,x)
=\sigma_{j,i}\,\frac{\partial u(t,x)}{\partial x_j},
\]
which means that a Neumann boundary condition
in the co-normal direction holds in (\ref{eq-pde-4})
on $\sigma S_{i,\eta}$ if $(i,\eta)\in\mathfrak{R}$.

We can thus solve (\ref{eq-pde-4}) by
reducing the problem to (\ref{eq-pde-3})
and use a Monte Carlo method in order to compute the values
of $u(t,x)$.

\subsection{The hypotheses}

Let us consider a domain $Q$ in $\mathbb{R}_+\times\mathbb{R}^d$.
For the sake of simplicity, we assume that $Q$
is the cylinder
$[0,T]\times D$ (with possibly $T=+\infty$),
where $D$ is an open, bounded domain of $\mathbb{R}^d$ with piecewise smooth
boundary.
Let us consider a function $a$ with values
in the space of $d\times d$-symmetric matrices
which is continuous on $D$
and everywhere positive definite, together with some functions
$b\dvtx Q\to\mathbb{R}^d$, $c\dvtx Q\to\mathbb{R}$ and $f\dvtx Q\to
\mathbb{R}$.
For all $(t,x)\in Q$, we denote by $\sigma(t,x)$ a $d\times
d$-symmetric matrix
such that $\sigma(t,x){\sigma}^{*}(t,x)=a(t,x)$.

We set
\[
L=\frac{1}{2}\sum_{i,j=1}^d a_{i,j}(t,x)\,\frac{\partial
^2}{\partial
x_i\,\partial x_j}
+\sum_{i=1}^d b_i(t,x)\,\frac{\partial}{\partial x_i}.
\]

Let us introduce the hypotheses needed to ensure
the convergence of our algorithm.
To set up a Monte Carlo numerical scheme,
one needs three inter-connected ingredients:
\begin{longlist}
\item The existence and the uniqueness of a
solution $u$ to the following PDE
%
%
\begin{equation}
\label{eq-pde-1}
\cases{
\dfrac{\partial u(t,x)}{\partial t}+Lu(t,x)\cr
\qquad{}+c(t,x)u(t,x)+f(t,x)=0, &\quad
on $[0,T]\times D$,\cr
u(T,x)=g(x), &\quad$x\in D$,\cr
u(t,x)=\phi(t,x), &\quad on $\Gamma_{\mathrm{d}}\subset[0,T)\times
\partial D$,\cr
\partial_n u(t,x)=0, &\quad on $\Gamma_{\mathrm{n}}\subset
[0,T)\times
\partial D$,}
\end{equation}
where $\partial_n$ denotes the co-normal derivative
along the lateral surface. $\Gamma_{\mathrm{d}}$ (respectively,
$\Gamma_{\mathrm{n}}$)
are subsets of $[0,T)\times\partial D$ on which
a Dirichlet (respectively, Neumann) boundary condition holds.

\item The existence of a solution to the diffusion
process associated with $L$. Note that since the simulation
involves distributions and not stochastic integrals,
we do not need strong
existence for the associated SDE.

\item The solution $u$ can be expressed
in terms of the probabilistic representation
%
%
\begin{eqnarray}
\label{eq-10}
u(t,x)&=&\mathbb{E}_{t,x}\biggl[\exp\biggl(\int_t^{\tau}
c(s,X_s)\,\vd s \biggr)\phi
(\tau,X_{\tau})\mathbh{1}_{\tau<T} \biggr]\nonumber\\
&&{} +\mathbb{E}_{t,x}\biggl[\exp\biggl(\int_t^{T} c(s,X_s)\,\vd s
\biggr)g(X_{T})\mathbh{1}_{\tau>T} \biggr]\\
&&{} +\mathbb{E}_{t,x}\biggl[\int_t^{\tau\wedge T} \exp\biggl
(\int_t^{s}
c(r,X_r)\,\vd r \biggr)f(s,X_s)\,\vd s \biggr],\nonumber
\end{eqnarray}
where $\tau$ is the first exit time
from $[0,+\infty)\times D$ by a point
of $\Gamma_{\mathrm{d}}$.
\end{longlist}
\begin{notation}
We denote by $\mathcal{P}$ the set of time--space parallelepipeds $P$
such that there exist $0\leq s<t\leq T$, $L_1,\ldots,L_d$ and
$x\in\mathbb{R}^d$ such that
\[
P=[s,t]\times\bigl(x+\widehat{\sigma}([-L_1,L_1]\times\cdots
\times[-L_d,L_d])\bigr),
\]
where $\widehat{\sigma}$ is a $d\times d$-matrix.
Possibly $t=+\infty$ (if $T=+\infty$).
\end{notation}

The assumptions that have to be done are the following:
\begin{enumerate}[(H1)]
\item[(H1)] There exists a subset $\mathcal{P}_D$ of $\mathcal{P}$
such that $Q=\bigcup_{P\in\mathcal{P}_D} P$.
Besides, if $P=[s,t]\times U\in\mathcal{P}$ for a parallelepiped $U$,
then for all $r\in[s,t)$, $[r,t]\times U\in\mathcal{P}$.
In other words, one can truncate the parallelepipeds in time.

\item[(H2)] There exist $\Gamma_{\mathrm{n}}$, $\Gamma_{\mathrm
{d}}$ contained in
$\partial Q=[0,T]\times\partial D$
and some subsets $\mathcal{P}_\mathrm{n}$, $\mathcal{P}_\mathrm{d}$
of $I$ such that
$\Gamma_{\mathrm{n}}\subset\bigcup_{P\in\mathcal{P}_\mathrm{n}}
\partial P$,
$\Gamma_{\mathrm{d}}\subset\bigcup_{P\in\mathcal{P}_\mathrm{d}}
\partial P$.
The closure of
$\Gamma_{\mathrm{n}}\cup\Gamma_{\mathrm{d}}$ is equal to
$[0,T]\times\partial D$
and $\Gamma_{\mathrm{n}}\cap\Gamma_{\mathrm{d}}=\varnothing$.
This means that the boundary of $[0,T]\times\partial D$
is split in two distinct parts,
where either the Dirichlet or the Neumann boundary
conditions hold. More precisely
a side of a parallelepiped in $\mathcal{P}_D$
contained in $\partial Q$ is either from $\Gamma_{\mathrm{n}}$ or
from $\Gamma_{\mathrm{d}}$.

\item[(H3)] The differential operator
$L$ is the generator of a continuous diffusion
process $X$ that is reflected at $\Gamma_{\mathrm{n}}$
and killed when hitting $\Gamma_{\mathrm{d}}\cup\{T\}\times D$.
The probabilistic representation of the solution
given by (\ref{eq-10}) holds
(see, e.g., \cite{lions84a} for existence results
of such reflected process
and \cite{stroock79a} if there are no reflections).

\item[(H4)] There exists an unique solution
$u$ of class $\mathcal{C}^{1,2}$ on $[0,T)\times D$
to (\ref{eq-pde-1}) which is continuous
on $[0,T]\times\overline{D}$.

\item[(H5)]
For a right parallelepiped $R$ and a matrix $\widehat{\sigma}$
let $P=[s,t]\times(x+\widehat{\sigma} R)\in\mathcal{P}_D$.
We associate with $P$ a vector $\widehat{b}\in\mathbb{R}^d$,
two constants $\widehat{c}$, $\widehat{f}$ and
we construct the differential operator
\[
\widehat{L} =\frac{1}{2}\sum_{k,l=1}^d\widehat{a}_{k,\ell}\,
\frac{\partial^2}{\partial x_k\,\partial x_\ell}
+\sum_{k=1}^d\widehat{b}_k\,\frac{\partial}{\partial_{x_k}}
\]
with $\widehat{a}=\widehat{\sigma}{\widehat{\sigma}}^{*}$.

Fix $\delta>0$. We assume that the solution $u$
to (\ref{eq-pde-1}) satisfies,
for any $y$ in the interior of $x+\widehat{\sigma}R$,
\[
\mathbb{E}_{s,y}
\biggl|
\int_s^{\tilde{\tau}} e^{\widehat{c}(r-s)}
\biggl(
\frac{\partial u}{\partial t}
+\widehat{L}u+\widehat{c}u-\widehat{f} \biggr)(r,\widehat{X}_r)
\,\vd r \biggr|\leq\delta,
\]
where $\widehat{X}$ is
the diffusion process generated by $\widehat{L}$,
and ${\tilde{\tau}}$ is its first exit time from $P$.
\end{enumerate}
\begin{remark}
If $T=+\infty$ and the coefficients are time-homogeneous
and $\Gamma_{\mathrm{d}}=[0,\infty)\times\gamma_{\mathrm{d}}$,
$\Gamma_{\mathrm{n}}=[0,\infty)\times\gamma_{\mathrm{n}}$, then
$v(x)=u(0,x)$
is solution to the elliptic PDE
%
%
\begin{equation}
\label{eq-pde-2}
\cases{
Lv(x)+c(x)v(x)=f(x), &\quad on $D$,\cr
v(x)=\phi(x), &\quad on $\gamma_{\mathrm{d}}\subset\partial D$,\cr
\partial_n v(x)=0, &\quad on $\gamma_{\mathrm{n}}\subset\partial D$.}
\end{equation}
Thus, by solving the parabolic PDE (\ref{eq-pde-1}),
we may also solve the elliptic PDE (\ref{eq-pde-2}).
We will thus focus only on (\ref{eq-pde-1}).
\end{remark}
\begin{remark}
The result of the existence
of a stochastic process reflected on some
part of the boundary of $[0,T)\times D$
is deduced from the existence of a stochastic
process reflected on the lateral boundary $[0,T)\times D$
which is killed when it hits $\Gamma_{\mathrm{n}}$.
\end{remark}

\subsection{The algorithm and its weak error}

In order to simplify the notation, if $T<+\infty$,
we denote the final condition $g$ of (\ref{eq-pde-1})
by $\phi(T,x)$.

Given $(t,x)\in Q$, the solution $u(t,x)$ of (\ref{eq-pde-1})
is computed by the Feynman--Kac formula. For this,
we have to simulate the diffusion process $X$
up to its first exit time $\tau$ from $Q$. We suppose here
that the particle cannot exit by a part of boundary
where a Neumann boundary condition holds. Let $u$
be the solution of (\ref{eq-pde-1}). Let us introduce the
following notation:
\[
\mbox{for }s\geq t\qquad
\cases{
\displaystyle Y_s=1+\int_t^s c(r,X_r) Y_r\,\vd r
=\exp\biggl(\int_t^s c(r,X_r)\,\vd r \biggr),\cr
\displaystyle Z_s=\int_t^s f(r,X_r) Y_r\,\vd r.}
\]
Then $u(t,x)$ is given by
%
%
\begin{equation}
\label{eq-300}
u(t,x)=\mathbb{E}_{t,x}[\phi(\tau,X_{\tau})Y_\tau+Z_\tau].
\end{equation}
We construct now the algorithm that approximates (\ref{eq-300})
by a Monte Carlo method.
%
%
\begin{algorithm}
\label{algo-2}
Assume that
we start initially at the point $(t,x)\in Q$
and fix a number $N$ of particles.
\begin{enumerate}[(1)]
\item[(1)] For $i=1,\ldots,N$ do
\begin{enumerate}[(A)]
\item[(A)] Set $(\theta_0,\Xi_0,Y_0,Z_0,W_0)=(t,x,1,0,1)$ and $k=0$.
\item[(B)] Repeat:
\begin{enumerate}[(a)]
\item[(a)] Choose an element $P^{(k)}\in\mathcal{P}_D$
of the form $P^{(k)}=[\theta_k,s]\times U$, $U\subset\mathbb{R}^d$
such that $(\theta_k,\Xi_k)$ belongs to the basis of $P$
($s$ is possibly infinite if, for example, $T=+\infty$
and the coefficients are time-inhomogeneous).
On $P^{(k)}$, consider the differential operator $L^{(k)}$
as well $c^{(k)}$ and $f^{(k)}$ which approximate
$L$, $c$ and $f$ as in (H5).
\item[(b)] Draw a realization of a random variable $(\theta_{k+1},\Xi_{k+1})$
with values in $(\{s\}\times U)
\cup((\theta_k,s)\times\partial U)$ and compute
its associated weight $w_k$ as shown
in Sections \ref{sec-para} and \ref{sub-reduc-para}
by considering the exit time and position from the
parallelepiped $P^{(k)}$.
\item[(c)] Compute $W_k=W_{k-1}w_k$ and
\begin{eqnarray*}
Y_{k+1}&=&Y_{k}\exp\bigl(c^{(k)} (\theta_{k+1}-\theta_k)\bigr),\\
Z_{k+1}&=&Z_{k}+f^{(k)}\int_{\theta_k}^{\theta_{k+1}}
\exp\bigl(c^{(k)} s\bigr)\,\vd s.
\end{eqnarray*}
\item[(d)] If $\Xi_{k+1}\in\Gamma_{\mathrm{d}}$ or $\theta
_{k+1}=T$, then exit
from the loop.
\item[(e)] Increase $k$.
\end{enumerate}
\item[(C)] Set
$(\theta^{(i)},\Xi^{(i)},Y^{(i)},Z^{(i)},W^{(i)})
=(\theta_{k+1},\Xi_{k+1},Y_{k+1},Z_{k+1},W_k)$.
\end{enumerate}
\item[(2)] Return
%
%
\begin{equation}
\label{eq-311}
\widehat{u}(t,x)=\frac{1}{N}\sum_{i=1}^N
\bigl(
W^{(i)}\phi\bigl(\theta^{(i)},\Xi^{(i)}\bigr
)Y^{(i)}+W^{(i)}Z^{(i)} \bigr).
\end{equation}
\end{enumerate}
\end{algorithm}


We denote from now on by $\widehat{\mathbb{P}}_x$ the distribution
of the Markov chain $\Lambda_k=(\theta_k,\Xi_k), k\geq0$.
Note that $(Y_k,Z_k,w_k)_{k\geq0}$ is obtained
from $(\Lambda_k)_{k\geq0}$.
\begin{proposition}
\label{prop-3}
For any $(t,x)\in[0,T)\times D$,
%
%
\begin{equation}
\label{eq-310}
|u(t,x)-\widehat{\mathbb{E}}_x[\widehat{u}(t,x)]|\leq
\delta\widehat{\mathbb{E}}_x[W_\nu\nu\exp(M\theta_\nu) ],
\end{equation}
where $\delta$ is defined in \textup{(H5)},
$\nu$ is the number of steps that $(\Lambda_k)_{k\geq0}$
takes to reach the boundary $\Gamma_{\mathrm{d}}\cap\{T\}\times D$
and
\[
M=\sup_{(s,y)\in[t,T)\times D} c(s,y).
\]
\end{proposition}
\begin{remark}
Note that the weak error in (\ref{eq-310}) does not
depend on the choice of the importance sampling technique
while the Monte Carlo error depends on this choice.
If the coefficients $a$, $b$, $f$ and $c$
are constant on the domain, one can choose $\delta=0$
and the simulation becomes exact.
\end{remark}
\begin{pf}
To the
Markov chain $(\Lambda_k)_{k\geq0}$ is associated a random sequence
of parallelepipeds $(P^{(k)})_{k=0,\ldots,\nu}$.
Let us denote by $\tau^{(k)}$ the successive times
the diffusion process $X$ reaches the boundary of the $P^{(k)}$'s.

Since $Z_0=0$, $Y_0=1$ and $u=\phi$ on the boundary of $Q$, we get
%
%
\begin{eqnarray}
\label{eq-302}\qquad
\widehat{\mathbb{E}}_x[\widehat{u}(t,x)]
&=& \widehat{\mathbb{E}}_x[W_\nu Y_\nu\phi(\theta_\nu,\Xi_{\nu})
+W_\nu Z_\nu] \nonumber\\
&=&u(t,x)+
\widehat{\mathbb{E}}_x\Biggl[W_\nu\sum_{k=0}^{\nu-1}
\bigl(Z_{k+1}-Z_k+Y_{k+1}u(\theta_{k+1},\Xi_{k+1})\\
&&\hspace*{177.2pt}{}-Y_{k}u(\theta_{k},\Xi
_{k})\bigr) \Biggr].\nonumber
\end{eqnarray}

Let $(\mathcal{G}_k)_{k\geq0}$ be the filtration generated by
the Markov chain $(\Lambda_k)_{k\geq0}$.
We remark that $Y_k$ and $Z_k$ are measurable with respect to
$\mathcal{G}_{k}$ while $w_k$ is measurable with respect to
$\mathcal{G}_{k+1}$ (since it is obtained
from $\theta_k$, $\Xi_k$, $\theta_{k+1}$ and $\Xi_{k+1}$).

By using the Markov property, after setting $W_{k+1,\nu}
=\widehat{\mathbb{E}}_x[w_{k+1}\cdots w_{\nu} | \mathcal{G}_{k+1}]$,
we get
\begin{eqnarray*}
&&\widehat{\mathbb{E}}_x[W_{\nu}(Z_{k+1}-Z_{k})]\\
&&\qquad=\widehat{\mathbb{E}}_x\bigl[W_{k+1,\nu} \widehat{\mathbb
{E}}_x[w_k
(Z_{k+1}-Z_{k}) | \mathcal{G}_{k}]W_{k-1}\bigr], \\
&&\widehat{\mathbb{E}}_x\bigl[W_\nu\bigl(Y_{k+1}u(\theta
_{k+1},\Xi_{k+1})-
Y_{k}u(\theta_{k},\Xi_{k})\bigr) \bigr] \\
&&\qquad=\widehat{\mathbb{E}}_x\bigl[W_{k+1,\nu} \widehat{\mathbb
{E}}_x\bigl[w_k
\bigl(Y_{k+1} u(\theta_{k+1},\Xi_{k+1})-Y_ku(\theta_{k},\Xi
_{k})\bigr) |
\mathcal{G}_{k}\bigr]W_{k-1}\bigr].
\end{eqnarray*}

Let us denote by $(X^{(k)},\mathbb{P}^{(k)}_{t,x})$
the process generated by the operator
$L^{(k)}$ with constant coefficients
$a^{(k)}$ and $b^{(k)}$
on $P^{(k)}$.
Define recursively
$(t^{(0)},x^{(0)})=(t,x)$ and
$(t^{(k+1)},x^{(k+1)})=({\overline{\tau}}^{(k)},X^{(k)}_{{\overline
{\tau}}^{(k)}})$
where ${\overline{\tau}}^{(k)}$ is, as above, the first exit time
from $P^{(k)}$
for the diffusion $X^{(k)}$.
Let also $f^{(k)}$ and $c^{(k)}$ be the values
that approach $f$ and $c$ on $P^{(k)}$,
and define also recursively $y^{(0)}=1$
and $y^{(k)}=y^{(k-1)}\exp(c^{(k)}(t^{(k+1)}-t^{(k)}))$.

By using the properties of $\widehat{\mathbb{P}}_x$
and the It{\^o} formula we obtain
\begin{eqnarray*}
&&\widehat{\mathbb{E}}_x\bigl[w_k\bigl(Y_{k+1}u(\theta_{k+1},\Xi
_{k+1})-Y_k
u(\theta_{k},\Xi_{k})\bigr) | \mathcal{G}_{k}\bigr] \\
&&\qquad=y^{(k)}\mathbb{E}^{(k)}_{t^{(k)},x^{(k)}}
\bigl[\bigl(e^{c^{(k)}(t^{(k+1)}-t^{(k)})}u\bigl
(t^{(k+1)},X^{(k+1)}_{t^{(k+1)}}\bigr)
-u\bigl(t^{(k)},x^{(k)}\bigr)\bigr)\bigr] \\
&&\qquad=y^{(k)}\mathbb{E}^{(k)}_{t^{(k)},x^{(k)}}
\biggl[\int_{t^{(k)}}^{t^{(k+1)}}e^{c^{(k)}(s-t^{(k)})} \biggl(\frac
{\partial
}{\partial t}+L^{(k)}+c^{(k)} \biggr)u\bigl(s,X_s^{(k)}\bigr)\,\vd s
\biggr].
\end{eqnarray*}

Also,
\[
\widehat{\mathbb{E}}_x[w_k (Z_{k+1}-Z_{k}) | \mathcal{G}_{k}]
=y^{(k)}\mathbb{E}^{(k)}_{t^{(k)},x^{(k)}}
\biggl[f^{(k)}\int_{t^{(k)}}^{t^{(k+1)}}e^{c^{(k)} s}\,\vd s \biggr].
\]
Under the hypothesis on the coefficients and the parallelepiped
$P^{(k)}$ we have
\begin{eqnarray*}
&&\bigl|\widehat{\mathbb{E}}_x\bigl[w_k\bigl(Y_{k+1}u(\theta
_{k+1},\Xi_{k+1})-Y_k
u(\theta_{k},\Xi_{k}) +Z_{k+1}-Z_k\bigr) | \mathcal{G}_{k}\bigr]
\bigr| \\
&&\qquad=
\biggl|y^{(k)}
\mathbb{E}^{(k)}_{t^{(k)},x^{(k)}}
\biggl[\int_{t^{(k)}}^{t^{(k+1)}}e^{c^{(k)}(s-t^{(k)})}\\
&&\hspace*{124pt}{}\times \biggl(
\biggl(\frac{\partial
}{\partial t}+L^{(k)} + c^{(k)} \biggr)u\bigl(s,X_s^{(k)}\bigr
)+f^{(k)} \biggr)\,\vd s \biggr] \biggr| \\
&&\qquad\leq y^{(k)}\delta
\leq
\widehat{\mathbb{E}}_x
[\delta w_k Y_{k} | \mathcal{G}_k ],
\end{eqnarray*}
since the $Y_k$'s (and so the $y^{(k)}$'s) are positive.
Hence, from (\ref{eq-302})
and the Jensen inequality applied to $|\cdot|$, we obtain
\[
|\widehat{\mathbb{E}}_x[Y_\nu\phi(\theta_\nu,\Xi_{\nu})+Z_\nu]
-\widehat{\mathbb{E}}_x[\widehat{u}(t,x)]|
\leq\delta
\widehat{\mathbb{E}}_x\Biggl[W_\nu\sum_{k=0}^{\nu-1} Y_k \Biggr].
\]
As $0<Y_k\leq e^{M\theta_k}$ for $k=0,\ldots,\nu$,
we deduce (\ref{eq-310}).
\end{pf}


\subsection{The Monte Carlo error}
\label{sub-monte-carlo-error}

In order to compute the solution $u(t,x)$ of (\ref{eq-pde-1}),
we have constructed the estimator
$\widehat{u}(t,x)$ given by (\ref{eq-311}) whose variance is
\[
\operatorname{Var}_{\widehat{\mathbb{P}}_x}\widehat{u}(t,x)
=\frac{1}{N}\operatorname{Var}_{\widehat{\mathbb{P}}_x}\bigl(W_\nu
\phi
(\theta_\nu,\Xi
_\nu)Y_\nu
+W_\nu Z_\nu\bigr).
\]
The Monte Carlo error depends on this
variance $s^2=\operatorname{Var}_{\widehat{\mathbb{P}}_x}\widehat{u}(t,x)$,
since asymptotically for $N\to\infty$
the true mean $\widehat{\mathbb{E}}_x[\widehat{u}(t,x)]$
lies in the interval
$[\widehat{u}(t,x)-2s,\widehat{u}(t,x)+2s]$ with
a confidence of $95.4 \%$.

We denote by $\widehat{\mathbb{P}}^{\mathrm{n}}$
the distribution
of $(\Lambda_k)_{k\geq0}$ with respect to the real distribution
of the exit time and position of the rectangles.
In this case the weights are equal to~$1$.
Any event $\Phi$ measurable with respect to $(\Lambda_k)_{k\geq0}$
satisfies $\widehat{\mathbb{P}}^{\mathrm{n}}[\Phi]=\widehat
{\mathbb{P}}_x[W\Phi]$.

We get thus
\[
\operatorname{Var}_{\widehat{\mathbb{P}}_x}\bigl(W_\nu\phi(\theta
_\nu
,\Xi_\nu)Y_\nu
+W_\nu Z_\nu\bigr)
=\Psi+\operatorname{Var}_{\widehat{\mathbb{P}}^{\mathrm
{n}}}(\widehat{u}(t,x))
\]
with
\[
\Psi=
\widehat{\mathbb{E}}^{\mathrm{n}}\bigl[(W_\nu-1) \bigl(\phi
(\theta_\nu,\Xi
_\nu)Y_\nu+Z_\nu\bigr)^2\bigr].
\]
This shows that a good choice for the density of the exit
time and position from the parallelepipeds is such that
$\Psi\leq0$ is as small as possible.
By the way, reducing the variance is a difficult task and
requires some automatic selection/optimization
techniques, as explained in the \hyperref[sec1]{Introduction}.

In addition, the numerical experiments we performed
up to now highlight another difficulty. $W_\nu$ may take large values,
and this implies meaningless values for $\widehat{u}(t,x)$.
That is why we suggest to keep track also of the
empirical distribution, or at least of the variance
of $W_\nu$.

In order to illustrate this, let us assume that the diffusion
process $X$ has no drift and that for the simulation,
the right parallelepipeds we use are squares centered
on the particle, and consider the same density
for the exit time and position. By a scaling argument,
the distribution of the weight $w_k$ at the $k$th step
does not depend on the size of the squares,
so that the $w_k$'s are independent and identically distributed
under $\widehat{\mathbb{P}}_x$.

Let us fix an integer $n$ such that $\nu\geq n$
a.s. (for example, the minimal number
of steps needed to reach an absorbing boundary).
We set $\xi^i=\log(w_i)$, so that
$W_n=\exp(\sum_{i=1}^n \xi^i )$.
As the $\xi^i$ are independent and identically distributed,
let us note $S_n=\sum_{i=1}^n \xi^i$, then
$S_n/\sqrt{n}$ converges to some normal random variable
$\chi$ with mean $m$ and variance $s^2$.
For $n$ large enough,
the distribution of $W_n$ is close to the distribution of
$\exp(\sqrt{n}\chi)$. We obtain,
with the expression of the Laplace transform for the normal
distribution, for $j\in\{1,2\}$,
\[
\widehat{\mathbb{E}}_x[(W_n)^j]\approx\mathbb{E}_x\bigl[\exp\bigl
(j\sqrt
{n}\chi\bigr)\bigr]
=\exp\biggl(m j\sqrt{n}+n\frac{j^2s^2}{2} \biggr).
\]
This leads us to the following approximation:
\begin{eqnarray*}
\operatorname{Var}_{\widehat{\mathbb{P}}_x}(W_n)
&\approx&\exp\bigl(2m\sqrt{n}+2ns^2 \bigr)
-\exp\biggl(m\sqrt{n}+\frac{n}{2}s^2 \biggr)\\
&\approx&
\exp(2ns^2)
\biggl(
\exp\biggl(\frac{m}{\sqrt{n}}+1 \biggr)
-\exp\biggl(\frac{m}{2\sqrt{n}}-\frac{3n}{2}s^2 \biggr)
\biggr)\\
&\mathop{\sim}\limits_{n\to\infty}& \exp(1+2ns^2).
\end{eqnarray*}
So, for large $n$, the variance
of $W_n$ explodes, while $\widehat{\mathbb{E}}_x[W_n]=1$
for any $n\geq1$.

In \cite{glynn} (see also \cite{glasserman04a}), Glynn and Iglehart exhibit
another argument that shows that the simulation performs
badly if too many steps are used.

\subsection{Population Monte Carlo}
\label{sec-population}

In order to overcome the explosion of the variance due to the weights
one can use a population Monte Carlo method.
This kind of method, also known as
quantum Monte Carlo,
sequential Monte Carlo, Green Monte Carlo$,\ldots$
has been used for a long time in physical simulations (see, e.g.,
\cite{iba} for a brief survey) but also in signal theory,
statistics$,\ldots.$
A probabilistic point of view is developed in the book \cite{delmoral04a}
of Del Moral.

In our case, instead
of simulating the particles one after another,
the idea is to keep track\vspace*{1pt} of the whole population
of $N$ particles $(y^{(i)})_{i\in\{1, \ldots, N\}}$
with time and space coordinates $(t^{(i)},x^{(i)})$ and a weight $w^{(i)}$
according to the algorithm given below. Each particle has two
possible states: \textit{still running} or \textit{stopped}. A
particle is stopped either at the first time it reaches
an absorbing boundary, or if its time is equal
to the finite final time $T$. Otherwise, the particle
is still running.
\begin{algorithm}
\label{algo-3}
This algorithm computes an approximation of the quantity $\mathbb
{E}_x[f(T\wedge\tau,X_{T\wedge\tau})]$
when $X_0=0$ by using a population of $N$ particles.

\begin{enumerate}
\item Set $n=0$; $n$ is the number of steps.
\item For $i$ from $1$ to $N$ set
\begin{enumerate}[(a)]
\item[(a)] $(w^{(i)}_0,t^{(i)}_0,x^{(i)}_0)=(0,0,x)$.
\end{enumerate}
\item Set $\mathfrak{S}=\varnothing$ and $\mathfrak{R}_n=\{
(w^{(i)}_0,t^{(i)}_0,x^{(i)}_0)\}_{i=1,\ldots,N}$.
\item While the set $\mathfrak{R}_n$ of still running particles
at step $n$ is nonempty do:
\begin{enumerate}[(a)]
\item[(a)] Set $\mathfrak{R}_{n+1}=\varnothing$.
\item[(b)] Do $\#\mathfrak{R}_n$ times the following operations:
\begin{enumerate}[(iii)]
\item[(i)] Pick a still running particle of index $j$ at random according
to a family of discrete probability distribution
\[
p_j=\frac{w^{(j)}_n}{
\sum_{k\ \mathrm{index}\ \mathrm{of}\ \mathrm{particles}\ \mathrm
{in}\ \mathfrak{R}_n} w^{(k)}_n},
\]
where $w^{(j)}_n$ is the weight of the particle
after $n$ iterations.
\item[(ii)] The particle is moved in time and space according
to the exit time and position from a time--space parallelepiped
that contains $(t^{(j)}_n,x^{(j)}_n)$.
Its new position is denoted $(t^{(j)}_{n+1},x^{(j)}_{n+1})$
and its associated weight $w^{(j)}_{n+1}$.
\item[(iii)] If $t^{(j)}_{n+1}=T$ or if $x^{(j)}_{n+1}$
belongs to an absorbing boundary, then
$(w^{(j)}_{n+1},t^{(j)}_{n+1},x^{(j)}_{n+1})$
is added to the set $\mathfrak{S}$ of stopped
particles. Otherwise, it is added to $\mathfrak{R}_{n+1}$.
\end{enumerate}
\item[(c)] Increment $n$ by $1$.
\end{enumerate}
\item Return
\[
\frac{1}{\sum_{i=1}^{N'} w^{(i)}}\sum_{i=1}^{N'}
w^{(i)}f\bigl(t^{(i)},x^{(i)}\bigr),
\]
when $\mathfrak{S}=\{(w^{(i)},t^{(i)},x^{(i)})\}_{i=1,\ldots,N'}$.
\end{enumerate}
\end{algorithm}

As we need to keep track of the positions of all the particles,
this algorithm is memory consuming. On the other hand,
it avoids the multiplication of the weights.
In addition, this algorithm can be modified in
the following way: instead of using $\#\mathfrak{R}_n$
particles at step $n$, it is possible to use $N$
particles, and in this case, one has to keep track
of the number of still running particles and to multiply
the weights by the proportion of still running particles.
The algorithm stops when the proportion of
still running particles is smaller than a given
threshold.
This approach can be used, for example,
for long time simulation, or to estimate
rare events, as, for example, in
\cite{cerou06a,delmoral04a,delmoral05a,lejay08a}.


\subsection{Estimation of the number of steps}
\label{sub-steps}

Let us consider now the estimation of the number
of steps. In order to do this we will use the techniques
employed in \cite{milstein95a,milstein99a,milstein04a}.

In Algorithm \ref{algo-2}, we have constructed the
Markov chain $(\Lambda_k)_{k\geq0}$
which is absorbed when reaching $\Gamma_{\mathrm{k}}=\Gamma_{\mathrm
{d}}\cap\{T\}\times D$.

For a function $u$ on $D$, we set
\[
Pu(t,x)=\widehat{\mathbb{E}}^{\mathrm{n}}[u(\Lambda_1) | \Lambda
_0=(t,x)]\quad\mbox{and}\quad A
=Pu(t,x)-u(t,x).
\]
The operator $A$ is the generator of a Markov chain.
\begin{lemma}
If $T<+\infty$ and
\[
\widehat{\mathbb{E}}^{\mathrm{n}}[\theta_1 | (\theta_0,\Xi
_0)=(t,x)]-t\geq\gamma,
\]
then
\[
\widehat{\mathbb{E}}^{\mathrm{n}}[\nu| (\theta_0,\Xi
_0)=(t,x)]\leq1+\frac{T-t}{\gamma}.
\]
\end{lemma}
\begin{pf}
Consider the problem
\[
\cases{
Av(t,x)=-g(t,x), &\quad on $Q$,\cr
u(t,x)=0, &\quad on $[0,T]\times\Gamma$,}
\]
whose solution is
\[
u(t,x)=\widehat{\mathbb{E}}^{\mathrm{n}}\Biggl[\sum_{k=0}^{\nu-1}
g(\Lambda_k) \Biggr].
\]
We remark that if $u$ and $g$ are well chosen this equality
gives a good estimate of $\widehat{\mathbb{E}}^{\mathrm{n}}[\nu]$.

Let $V(t,x)$ be the function $V(t,x)=(T-t)\mathbh{1}_{(t,x)\in Q}$.
For $(t,x)$ in $Q$, we have
\[
AV(t,x)=\widehat{\mathbb{E}}^{\mathrm{n}}[V(\theta_{1},\Xi_1) |
(\theta_0,\Xi_0)=(t,x)]
-(T-t)\leq-\gamma.
\]
Hence $T-t\geq
\widehat{\mathbb{E}}^{\mathrm{n}}[\sum_{k=0}^{\nu-1} \gamma|
(\theta_0,\Xi_0)=(t,x)]$
and the result follows easily.
\end{pf}
\begin{lemma}
With the previous notation, for every $L>0$ fixed, we have
\[
\sup_{x\in Q}\widehat{\mathbb{P}}^{\mathrm{n}}[\nu\geq L |
(\theta_0,\Xi_0)=(t,x) ]
\leq(1+T-t)\exp\bigl(-c\gamma L/(1+T-t)\bigr),
\]
where $c$ is a constant depending on $\gamma$; more precisely $c$
converges to $1$ as $\gamma$ decreases to $0$.
\end{lemma}
\begin{pf}
The proof follows from the one of Theorem 7.2 in \cite{milstein99a}.
\end{pf}
\begin{lemma}
If $T=+\infty$, $Q$ is bounded and
\[
\widehat{\mathbb{E}}^{\mathrm{n}}[|x+\Xi_1+c|^2]\geq\gamma>0,
\]
where $c$ is such that
$\min_{x\in\overline{Q}}|x+c|\geq C>0$. Then
\[
\widehat{\mathbb{E}}^{\mathrm{n}}[\nu]\leq\frac
{B^2-C^2}{B^2-\gamma}
\]
with $B>\max\{\gamma,\sup_{x\in\overline{Q}}|x+c|\}$.
\end{lemma}
\begin{pf}
Let us proceed as in \cite{milstein95a}.
Choose a vector $c$ such that $\min_{x\in\overline{Q}}|x+c|\geq C>0$,
and set
\[
V(t,x)=\cases{
B^2-|x+c|^2, &\quad if $(t,x)\in\mathbb{R}_+\times Q$,\cr
0, &\quad otherwise.}
\]
Thus for $B^2 > \gamma$,
\[
AV(t,x)\leq
|x+c|^{2}-\widehat{\mathbb{E}}^{\mathrm{n}}[|x+\Xi_1+c|^2 | (\theta
_0,\Xi_0)=(t,x)]
\leq B^2-\gamma
\]
and the result follows.
\end{pf}


\section{Numerical examples}
\label{sec-numerical}

We present in this section some numerical examples
in order to test our algorithm.


\subsection{Speeding up the random walk on squares algorithm}

In \cite{milstein99a} (see also \cite{milstein04a}), Milstein
and Tretyakov
propose a method to simulate Brownian motions
and solutions of SDEs by using the first exit
time and position from a hyper-cube
or a time--space parallelepiped with cubic space basis.
A similar method has been previously proposed by Faure
in his Ph.D. thesis \cite{faure92a}.
This method is a variation of the random walk
on spheres method.
Some authors already used random walk
on squares and rectangles by using the explicit expression
of the Green function but without simulating
the exit time (see, e.g., \cite{simonov04a}).
One of the main features of our approach
is the simulation of the couple of nonindependent
random variables (exit time, exit position)
by means of real valued random variables.
We have
explained in \cite{deaconulejay05a} how to extend this
approach to rectangles and the starting point everywhere
in the rectangle. This approach is still using only
one-dimensional distributions.
However, by using symmetry properties, we can notice that
it is simpler
to deal with squares centered on the current position
of the particle than with a rectangle.

Nevertheless, the computation may
be time consuming.
We are looking now to speed up the computations by
using a simple density for the exit position.

Let us consider here the $d$-dimensional hypercube $C=[-1,1]^d$,
and a fixed time $T>0$ (possibly $T=+\infty$).
Let $B$ be a $d$-dimensional Brownian motion.
We set $\tau^B=\inf\{t>0 | B_t\notin C\}$.
Let $W$ be a one-dimensional Brownian motion.
We set $\tau^W_{[-1,1]}=\inf\{t>0 | W_t\notin[-1,1]\}$,
$R(t)=\mathbb{P}_0[\tau^W_{[-1,1]}<t]$,
$r$ the density of $\tau^W_{[-1,1]}$,
$S(t,y)=\mathbb{P}_0[W_t<y | t<\tau^W_{[-1,1]}]$
and $s(t,y)=\partial_y S(t,y)$ the density
of $W_t$ given $\{t<\tau^W_{[-1,1]}\}$.

Let us note that we can easily switch from $C=[-1,1]^d$
to any hypercube $[-L,L]^d$ after a
scaling argument in space and time. Thus, from a numerical
point of view, we need only to implement
the required functions $r$, $s$, $R$ and $S$ on $[-1,1]$.
Analytical
expressions for these distribution functions are easily deduced
from the series presented in the \hyperref[app-1]{Appendix}.

To simulate the exit time and position
from $[0,T]\times C$, we proceed in
the following steps:
\begin{itemize}
\item Compute the probability
$\beta=1-(1-R(T))^d$ that $\tau^B<T$.
\item With probability $\beta$,
decide if $\{\tau^B<T\}$ happens or not.
\item If $\{\tau^B<T\}$ happens:
\begin{itemize}
\item For a realization $\mathsf{U}$ of a uniform
random variable $U$ on $[0,1)$, set
\[
\overline{\tau}^B=R^{-1}\bigl(1-(1-\mathsf{U}\beta)^{1/d}\bigr),
\]
which is a realization of $\tau^B$ given $\{\tau^B<T\}$.
\item Choose with probability $1/2d$
an exit side $(J,\varepsilon)$, and set $\xi_J=\varepsilon$.
\item For each $i=1,\ldots,d$, $i\not=J$,
set $\chi_{i}=\sqrt{\mathsf{U}_{i}}$,
where the $\mathsf{U}_{i}$'s are $d-1$ independent realizations
of uniform random variables on $[0,1)$.
With probability $1/2$, set $\xi_{i}
=\chi_{i}-1$ and with probability $1/2$,
set $\xi_{i}=1-\chi_{i}$.
\item Compute the weight
\[
w=\frac{1}{1-R(\overline{\tau}^B)}\prod_{i=1,\ldots,d,i\neq J}
\frac{s(\overline{\tau}^B,\xi_{i})}{\chi_i}.
\]
\end{itemize}
\item If $\{\tau^B\geq T\}$ happens, then:
\begin{itemize}
\item Set $\overline{\tau}^B=T$.
\item For $i=1,\ldots,d$,
set $\chi_{i}=\sqrt{\mathsf{U}_{i}}$,
where the $\mathsf{U}_{i}$'s are $d-1$ independent realizations
of uniform random variables on $[0,1)$.
With probability $1/2$, set $\xi_{i}
=\chi_{i}-1$ and with probability $1/2$,
set $\xi_{i}=1-\chi_{i}$.
\item Compute the weight
\[
w=\frac{1}{1-\beta}\prod_{i=1,\ldots,d}\frac{s(T,\xi_i)}{\chi_i}.
\]
\end{itemize}
\end{itemize}
$(\overline{\tau}^B,\xi_1,\ldots,\xi_d)$
represent the first exit time and position
from $[0,T]\times C$, and $w$ is the associated
weight.

For the random walk on squares we can also use
the idea proposed in \cite{milstein99a}
and in \cite{faure92a}. This leads to the following
algorithm:
\begin{itemize}
\item Compute the probability
$\beta=1-(1-R(T))^d$ that $\tau^B<T$.
\item With probability $\beta$,
decide if $\{\tau^B<T\}$ happens or not.
\item If $\{\tau^B<T\}$ happens:
\begin{itemize}
\item For a realization $\mathsf{U}$ of a uniform
random variable $U$ on $[0,1)$, set
\[
\overline{\tau}^B=R^{-1}\bigl(1-(1-\mathsf{U}\beta)^{1/d}\bigr),
\]
which is a realization of $\tau^B$ given $\{\tau^B<T\}$.
\item Choose with probability $1/2d$
an exit side $(J,\varepsilon)$, and set $\xi_J=\varepsilon$.
\item For each\vspace*{-3pt} $i=1,\ldots,d$, $i\not=J$,
draw $\xi_{i}$ according to the
distribution of $B^i_{\overline{\tau}^B}$ given $\tau
^{B^i}>\overline{\tau}^B$,
where $\tau^{B^i}=\inf\{t>0 | B^i\notin[-1,1]\}$.
\end{itemize}
\item If $\{\tau^B\geq T\}$ happens, then:
\begin{itemize}
\item Set $\overline{\tau}^B=T$.
\item For $i=1,\ldots,d$,
draw $\xi_{i}$ according to the
distribution of $B^i_{\overline{\tau}^B}$ given $\tau
^{B^i}>\overline{\tau}^B$.
\end{itemize}
\end{itemize}

In both cases, we use tabulated values for $R$ and $R^{-1}$.
In order to simulate $B^i_t$ given $\tau^{B^i}>t$,
we use the rejection method proposed by Faure
in \cite{faure92a} for $t\in[0.25,2]$. Otherwise,
we draw $B^i_t$ by using the fact that it is
equal to $S^{-1}(t,U)$ for some random variable $U$
with uniform distribution on $[0,1)$.
This is the method proposed by Milstein and
Tretyakov in \cite{milstein99a}.
For $t>2$, the latter method is more efficient
than the previous one. For $t<0.2$, the rejection method
may give wrong results. For $t$ close to $0.2$,
the rejection method can be up to $6$ times faster
than the inversion method, while for $t$ close to $2$,
they are comparable in the computation time.

If the Brownian motion reaches the side labeled by $(1,-1)$
first at time $\tau^B$, then in order to simulate
$B^i_t$ for $i=2,\ldots,d$ we use a random variable with density
$\phi(x)=1+x$ if $x\in(-1,0]$ and $\phi(x)=1-x$ if $x\in[0,-1)$.
In this case, the weights $w$ are close to $1$
as we see in Table \ref{table-4},
%
%
\begin{table}
\caption{Speeding up the random walk on squares:
experiments
with 1,000,000 samples are used}\label{table-4}
\begin{tabular*}{\tablewidth}{@{\extracolsep{\fill}}lcccd{2.2}@{}}
\hline
\textbf{Method} & $\bolds T$ & \textbf{Mean of} $\bolds w$ & \textbf
{Variance of} $\bolds w$
& \multicolumn{1}{c@{}}{\textbf{Time (s)}}\\
\hline
Walk on squares & $0.1$ & -- & -- & 94\\
Imp. sampling & $0.1$ & $1.0005$ & $0.28$\phantom{0} & 3.2\\
Walk on squares & $0.2$ & -- & -- & 82\\
Imp. sampling & $0.1$ & $0.9997$ & $0.014$ & 1.8\\
Walk on squares & $0.5$ & -- & -- & 10\\
Imp. sampling & $0.5$ & $0.9999$ & $0.021$ & 1.2\\
Walk on squares & $1.0$ & -- & -- & $10$\\
Imp. sampling & $1.0$ & $0.9994$ & $0.017$ & 1\\
Walk on squares & $+\infty$ & -- & -- & 10\\
Imp. sampling & $+\infty$ & $0.9998$ & $0.015$ & 0.98\\
\hline
\end{tabular*}
\end{table}
and the execution time is usually divided by $10$.
For $T=0.1$, the variance of $w$ is
too high and leads to some instabilities.
In this case, it is preferable
to simulate the exact distributions of $B_T$
given $\{T\leq\tau^B\}$.


\subsection{Solving a bi-harmonic problem}

To test the validity of our approach with respect
to other algorithms,
we consider first an example borrowed in \cite{milstein99a}
(see also \cite{milstein04a}, page 332).
Let $D=[-1,1]^2$, and consider the bi-harmonic equation
%
%
\begin{equation}
\label{eq-ex-num-bih-1}
\cases{
\frac{1}{2}\triangle^2 u(x)=1, &\quad$x\in D$,\cr
u(x)=\phi(x), &\quad on $\partial D$,\cr
\frac{1}{2}\triangle u(x)=\psi(x), &\quad on $\partial D$,}
\end{equation}
with
%
%
\begin{eqnarray}
\phi(x_1,\pm1)&=&\frac{1+x_1^4}{12},\qquad\phi(\pm1,x_2)=
\frac{1+x_2^4}{12},\\
\psi(x_1,\pm1)&=&\frac{1+x_1^2}{2},\qquad\psi(\pm1,x_2)=
\frac{1+x_2^2}{2}.
\end{eqnarray}
After setting $v(x)=\frac{1}{2}\triangle u(x)$, (\ref{eq-ex-num-bih-1})
may be transformed into the system
\[
\cases{
\frac{1}{2}\triangle v(x)=1\mbox{ on $D$}, &\quad
with $u(x)=\psi(x)\mbox{ on }\partial D$,\vspace*{2pt}\cr
\frac{1}{2}\triangle u(x)-v(x)=0\mbox{ on $D$},&\quad
with $u(x)=\phi(x)\mbox{ on }\partial D$,}
\]
whose exact solution is
\[
u(x)=\frac{x_1^4+x_2^4}{12},\qquad
v(x)=\frac{x_1^2+x_2^2}{2}.
\]
By It{\^o}'s formula, it is easy to show that
\begin{eqnarray*}
u(x)&=&\mathbb{E}[\phi(x+B_{\tau^B})]-\mathbb{E}[{\tau^B}\psi
(x+B_{\tau^B})]+\tfrac{1}{2}
\mathbb{E}[({\tau^B})^2],\\
v(x)&=&\mathbb{E}[\psi(x+B_{\tau^B})]-\mathbb{E}[{\tau^B}],
\end{eqnarray*}
where $B$ is a two-dimensional Brownian motion, and ${\tau^B}$ is, as
above, its first exit time from $D$.

Here, in contrast with the values presented in \cite{milstein99a},
we only need to use one square, since we are not forced to
start from its center. We compare the results given
by our algorithm (first lines) with the one given by the random walk
on rectangles (second line).
Each side is chosen uniformly with probability $1/4$.
The time is drawn by using an exponential random
variable of parameter $1/(1-\varepsilon x_i)$
if $(i,\varepsilon)$ is the exit side. The position
is drawn uniformly on the exit side. This strategy corresponds in some sense
to a ``naive'' and simple way to choose the
exit time and position.

As we evaluate quantities of the form $\mathbb{E}[f({\tau^B},B_{\tau^B})]$,
we report the quantities $\mu_n\pm2\sigma_n/\sqrt{n}$,
where $\mu_n$ is the empirical mean
of $f({\tau^B},B_{\tau^B})$ with $n$ samples,
and $\sigma_n$ is the corresponding empirical standard deviation.
The interval $[\mu_n-2\sigma_n/\sqrt{n},\mu_n+2\sigma_n/\sqrt{n}]$
represents the $95.5 \%$ confidence interval for $\mathbb{E}[f({\tau
^B},B_{\tau^B})]$.
The estimations $\overline{u}(x)$ and $\overline{v}(x)$
of $u$ and $v$ for three points are given in Table \ref{table-3}.

%
%
\begin{table}[b]
\caption{Solution of the bi-harmonic equation:
the first line of each row contains the results for our algorithm, and
the second line contains the results for the random walk on
rectangles}\label{table-3}
\begin{tabular*}{\tablewidth}{@{\extracolsep{\fill}}lcccccd{4.2}@{}}
\hline
$\bolds x$ & $\bolds n$ & $\bolds{u(x)}$ & $\bolds{\overline{u}(x)}$
& $\bolds{v(x)}$ & $\bolds{\overline{v}(x)}$
& \multicolumn{1}{c@{}}{\textbf{Time (s)}}\\
\hline
$(0.3,0.5)$ & $10^4$ & $0.00588$ & $0.0047\pm0.0037$ & $0.17000$ & $0.1638\pm0.0081$ & 0.03\\
& & & $0.0064\pm0.0039$ & & $0.1684\pm0.0081$ & 3.8\\
-- & $10^5$ & -- & $0.0061\pm0.0012$ & -- & $0.1669\pm0.0026$ & 0.23
\\
& & & $0.0062\pm0.0012$ & & $0.1679\pm0.0026$ & 38 \\
-- & $10^6$ & -- & $0.0059 \pm0.0004$ & -- &  $0.1698\pm0.0008$ &
2.2\\
& & & $0.0059\pm0.0004$ & & $0.1696\pm0.0008$ & 381\\
$(0.7,0.8)$ & $10^4$ & $0.05414$ & $0.0480\pm0.0017$ & $0.56500$ &
$0.5297\pm0.0064$ & 0.02 \\
& & & $0.0553 \pm0.0020$ & & $0.5707\pm0.0061$ & 7
\\
-- & $10^5$ & -- & $0.0526\pm0.0005$ & -- & $0.5593\pm0.0019$  &
0.2\\
& & & $0.0543\pm0.0006$ & & $0.5652 \pm0.0019$ & 73\\
-- & $10^6$ & -- & $0.0536\pm0.0002$ & -- & $0.5654\pm0.0006$
& 2.5\\
& & & $0.0542\pm0.0002$ & & $0.5650\pm0.0006$ & 726\\
$(0.9,0.9)$ & $10^4$ & $0.10935$ & $0.1103\pm0.0009$ & $0.81000$ &
$0.8186\pm0.0034$ & 0.01
\\
& & & $0.1109\pm0.0020$ & & $0.8105\pm0.0038$ & 11\\
-- & $10^5$ & -- & $0.1131\pm0.0002$ & -- &
$0.8390\pm0.0006$ & 0.2
\\
& & & $0.1095\pm0.0003$ & & $0.8107\pm0.0011$ & 112\\
-- & $10^6$ & -- & $0.1087\pm0.0001$ & -- &
$0.8097\pm0.0003$ & 2
\\
& & & $0.1093\pm0.0001$ & & $0.8100\pm0.0003$ & 1100\\
\hline
\end{tabular*}
\end{table}

Although a small numerical bias seems to appear,
our algorithm provides results comparable with
the random walk on rectangles method. The execution
time is much smaller than the one given by this method
(also the one given
by the random walk on squares, for which
the simulation of one step takes less time, but where more steps are needed).


\subsection{Estimation of rare events: Computing hitting probabilities}

Let us consider the following problem: what is the probability $p(x)$
that starting from a point $x$ in a domain $D$ a Brownian
motion reaches a part $S$ of the boundary $\partial D$?
It is well known that $p$ is the solution of the Dirichlet
problem
%
%
\begin{equation}
\label{eq-simul-rare-event-1}
\tfrac{1}{2}\triangle p(x)=0\mbox{ on }D
\quad\mbox{and}\quad
p(x)=\cases{
1, &\quad if $x\in S$,\cr
0, &\quad if $x\in\partial D\setminus S$.}
\end{equation}

\begin{figure}

\includegraphics{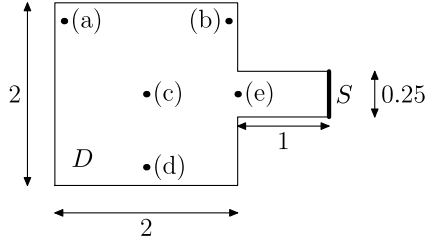}

\caption{A simple domain $D$.}\label{fig-2}
\end{figure}

We illustrate our method on the simple two-dimensional domain $D$
drawn in Figure \ref{fig-2} and we compute
the value of $p$ at the five points marked, respectively, by (a), (b),
(c), (d)
and (e) on Figure \ref{fig-2}.

%
%

\begin{figure}[b]

\includegraphics{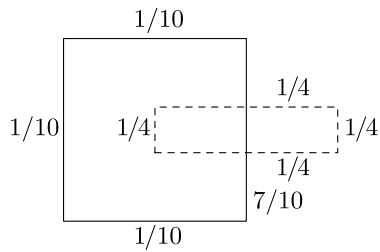}

\caption{Decomposition of $D$ into rectangles.}\label{fig-3}
\end{figure}

To set up our algorithm, we use two rectangles
as in Figure \ref{fig-3}. The numbers marked on each side
are the probabilities to reach each
one of these sides.

In order to obtain the simulated exit time we draw an exponential
random variable with parameter $\alpha$ where $\alpha$ is given
by $\alpha= 1/ (\sqrt{L_i/2})$. The $L_i$ notes the length of the
rectangle in the direction perpendicular to the boundary that the particle
hits.

%
%

We perform 100,000 samples; each
computation takes around 1 s on our computer
(a MacBook 12$''$, 2 GHz with a code written in \texttt{C}).
The values for $p$ are given in Table \ref{table-1}.
We perform a comparison with the value given
by MATLAB/PDEtool where
(\ref{eq-simul-rare-event-1}) is solved by using
a finite element method, and with the
method of random walk on rectangles \cite{deaconulejay05a}
which is exact (up to the Monte Carlo error), for such a domain.
In this case, with a sample of size $n$, the variance
of the empirical mean is $p(x)(1-p(x))/n$.

%
%
\begin{table}
\tablewidth=236pt
\caption{Computation of $p(x)$ at given points of $D$}\label{table-1}
\begin{tabular*}{\tablewidth}{@{\extracolsep{\fill}}lccc@{}}
\hline
\textbf{Point} & \textbf{Import. sampling} & \textbf{Finite element}
& \textbf{Walk on rect.}\\
\hline
(a) & $3.32\cdot10^{-6}$ & $3.39\cdot10^{-6}$ & $0.00$\phantom{$1\cdot0^{-6}$}\\
(b) & $2.31\cdot10^{-5}$ & $2.23\cdot10^{-5}$ & $1.00\cdot10^{-5}$\\
(c) & $1.70\cdot10^{-4}$ & $1.77\cdot10^{-4}$ & $1.90\cdot10^{-4}$\\
(d) & $4.43\cdot10^{-5}$ & $4.64\cdot10^{-5}$ & $3.00\cdot10^{-5}$\\
(e) & $2.79\cdot10^{-3}$ & $2.81\cdot10^{-3}$ & $2.36\cdot10^{-3}$\\
\hline
\end{tabular*}
\end{table}

We notice that the results given by our method are
close to the one given by
the finite element method. As one can expect, the random walk
on rectangles (and any other methods that do not
rely on importance sampling or variance reduction
techniques) is not efficient to estimate
the values of $p(x)$ when they are of the same order
as the standard deviation of the empirical mean.

\begin{figure}[b]

\includegraphics{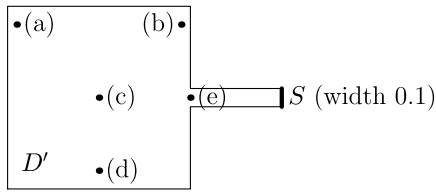}

\caption{A simple domain $D'$.}\label{fig-4}
\end{figure}

In order to test the validity of our method
for the simulation of rare events, we use
the domain $D'$ as in Figure \ref{fig-4}.

%
%

The numerical results are reported in Table \ref{table-2}.
$p_n$ is the empirical mean
with $n=100\mbox{,}000$ samples, and $s_{50}(p_n)$
is the empirical standard deviation computed
over $50$ realizations of $p_n$.
We obtain really good results even while computing
small probabilities
of order of magnitude $10^{-10}$.


\subsection{Simulation of SDEs: Approximation close to the boundary}

\begin{figure}[b]

\includegraphics{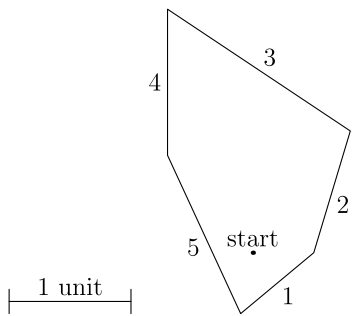}

\caption{Domain $D$ with the label of the sides and the starting
point.}\label{fig-5}
\end{figure}

Let us consider the two-dimensional SDE
\[
X_t=x+\int_0^t \sigma(X_s)\,\vd B_s
\qquad\mbox{with }\sigma(x)= \left[\matrix{
1 & \frac{1}{2}\sin(x_1+x_2)\vspace*{2pt}\cr
0 & 1} \right],
\]
which is driven by a two-dimensional Brownian motion $B$.
The process $X$ is killed when it exits from the
domain $D$ which is represented in Figure \ref{fig-5}.

In order to simulate $X$, we use either an Euler
scheme with a time step of $0.0025$ or a (possibly modified) random
walk on squares. The squares sides lengths are smaller
than $2L$ with $L=0.05$ (note that the time step of the
Euler scheme corresponds to $0.05^2$ which is close to
the average exit time of the square $[0.1,0.1]^2$).
As the diffusion moves in a bounded domain, we
use to deal with the boundary condition and apply the technique
proposed in \cite{lejaym}: if
the distance between the position of the particle
and the boundary is smaller than $2L$, we choose
the square such that one of its sides is included
in the boundary when it is possible to do so.

Unless the coefficients of the SDE are constant,
one needs to simulate many couples of
exit times and positions from small squares, and
the computational time becomes very large
and is not competitive with respect to the Euler
scheme. In addition, when the random walk on
squares is coupled with importance sampling,
the weights grow quickly (see Section \ref{sub-monte-carlo-error}).

When the Euler scheme is used, we simply stop
the algorithm when the particle leaves the domain $D$.
This is a crude way to proceed, and some refinements
can be done (see, e.g., \cite{gobet00a}).
Note that the exit time is then overestimated.

%
%
\begin{table}
\caption{Computation of $p(x)$ at given points of $D'$}\label{table-2}
\begin{tabular*}{\tablewidth}{@{\extracolsep{\fill}}lccc@{}}
\hline
\textbf{Point} & $\bolds{p_n}$ & $\bolds{s_{50}(p_n)}$ & \textbf{Finite element} \\
\hline
(a) & \phantom{$^{0}$}$1.00\cdot10^{-10}$ & $2.3\cdot10^{-11}$ & \phantom{$^{0}$}$1.15\cdot10^{-10}$
\\
(b) & \phantom{$^{0}$}$7.67\cdot10^{-10}$ & $1.6\cdot10^{-10}$ & \phantom{$^{0}$}$8.13\cdot10^{-10}$
\\
(c) & $5.19\cdot10^{-9}$ & $1.0\cdot10^{-9}$\phantom{$^{0}$} & $6.61\cdot10^{-9}$ \\
(d) & $1.31\cdot10^{-9}$ & $2.8\cdot10^{-10} $ & $1.73\cdot10^{-9}$
\\
(e) & $2.27\cdot10^{-7}$ & $4.9\cdot10^{-8}$\phantom{$^{0}$} & $2.29\cdot10^{-7}$ \\
\hline
\end{tabular*}
\end{table}

%
%

The idea is to mix the two methods and
to use the Euler scheme inside the domain,
and a random walk on squares when the particle
is close to the boundary. We improve thus the simulation
as in this case the behavior of the particle
is taken into account.
In addition, it is
possible by making a change of measure, to
increase or to decrease the probability that the
particle hits the boundary.

Our aim is here to increase the number of
particles which are not killed before a given
time $T$.
When one side of the square is set on the boundary,
we use a probability $p$ that the particle reaches
the side of the square that is opposite to the
boundary, and $q=(1-p)/3$ for any other
side. We have thus a ``repulsing'' effect.

We use $P_1=\{p=0.7,q=0.1\}$ and $P_2=\{p=0.91,q=0.03\}$.

In order to avoid the explosion of the variance
of the weight, we have used a limitation $N_{\max}$ for
the number of times this procedure is used. The variance
of the weight for each time this procedure is used
is $0.04$ for the set $P_1$
and $0.34$ for the set $P_2$.

%
%
\begin{table}
\tabcolsep=0pt
\caption{Simulations of the proportions (in \%) of
the particles
reaching a given part of the boundary as well as the surviving
particles at time $T$. We write ``unstable'' in the column for the variance
of weights when the mean of the global weights is rather far from
$1$}\label{table-6}
\begin{tabular*}{\tablewidth}{@{\extracolsep{\fill}}lcccccccccc@{}}
\hline
$\bolds T$ & \textbf{Type} & \textbf{Side 1} & \textbf{Side 2} & \textbf{Side 3}
& \textbf{Side 4} & \textbf{Side 5} & \textbf{Final time} &
\textbf{Var. weights} & $\bolds{N_{\max}}$ & \textbf{Time} \\
\hline
\multicolumn{11}{c}{{Test with set of probabilities $P_1$ on
the boundary}}\\
[4pt]
1 & Est. & $31.55$ & $12.39$ & $4.16$ & $0.44$ & $51.27$ & $0.17$ &
\phantom{00}$9.3$\phantom{0} & \phantom{00}$5$ & \phantom{0}$73$\\
& Sim. & $31.91$ & $12.97$ & $5.07$ & $0.65$ & $48.96$ & $0.41$ & & \\
1 & Est. & $30.81$ & $13.19$ & $4.06$ & $0.32$ & $51.43$ & $0.17$ &
\phantom{0}$29.9$\phantom{0} & \phantom{0}$10$ & \phantom{0}$83$\\
& Sim. & $31.93$ & $13.08$ & $5.42$ & $0.75$ & $48.10$ & $0.66$ & & \\
1 & Est. & $31.05$ & $13.83$ & $4.37$ & $0.42$ & $50.14$ & $0.17$ &
\phantom{0}$30.0$\phantom{0} & \phantom{0}$20$ & \phantom{0}$93$\\
& Sim. & $32.01$ & $13.40$ & $5.57$ & $0.91$ & $47.10$ & $0.96$ & & \\
1 & Est. & $30.96$ & $13.54$ & $4.08$ & $0.36$ & $50.84$ & $0.19$ &
\phantom{0}$56.55$ & $100$ & \phantom{0}$99$\\
& Sim. & $31.78$ & $13.27$ & $5.57$ & $0.98$ & $46.83$ & $1.36$ & & \\
[4pt]
\multicolumn{11}{c}{{Test with set of probabilities $P_2$ on
the boundary}}\\
[4pt]
1 & Est. & $29.45$ & $12.11$ & $3.49$ & $0.58$ & $54.19$ & $0.14$ &
$426$\phantom{00.} & \phantom{00}$5$ & \phantom{0}$90$\\
& Sim. & $32.13$ & $13.07$ & $5.71$ & $0.81$ & $47.61$ & $0.95$ & & \\
1 & Est. & $33.76$ & $11.78$ & $5.71$ & $0.37$ & $48.16$ & $0.19$ &
$65.5$ (unstable) & \phantom{0}$10$ & $117$\\
& Sim. & $32.03$ & $13.50$ & $6.70$ & $1.13$ & $45.21$ & $1.48$ & & \\
1 & Est. & $31.28$ & $14.19$ & $3.75$ & $0.44$ & $50.10$ & $0.21$ &
\phantom{0}$51.08$ (unstable) & \phantom{0}$20$ & $162$\\
& Sim. & $31.18$ & $13.48$ & $7.64$ & $1.62$ & $42.44$ & $3.62$ & & \\
1 & Est. & $29.87$ & $13.73$ & $2.83$ & $0.30$ & $53.03$ & $0.23$ &
$312.5$ (unstable)\phantom{0} & $100$ & $223$\\
& Sim. & $28.13$ & $12.21$ & $7.50$ & $1.58$ & $36.36$ & $14.23$\phantom{0} & & \\
\hline
\end{tabular*}
\end{table}

All the simulations are done with 100,000 particles.
The results are summarized in Table \ref{table-6}.
For $T=1$, the proportion of particles still alive
is of order $0.19\%$ (using the Euler scheme without specific
treatment on the boundary, we get an estimation of $0.33\%$,
yet for a quicker simulation of $7$ s).
With a population Monte Carlo method, we obtain an estimate
of $0.17\%$, using the set $P_1$ and a running time of $126$ s.
We see that our scheme allows
one to get much more alive particles.

%
\begin{appendix}\label{app-1}
\section*{Appendix: How to get densities for different situations?}

We present in this section analytical expressions for
the density in different cases.

Except for the case of a drifted Brownian motion with
Dirichlet boundary condition at one endpoint of $[-1,1]$
and a Neumann boundary condition at the other endpoint of $[-1,1]$,
we obtain two expressions, one which follows from
the images method and the other one from
the spectral decomposition. From a numerical point
of view, the spectral decomposition gives rise to series
that converge very quickly for large times. It is worth using the
expressions given by the
method of images for small times.

\subsection{Brownian motion without drift}

We are interested in this section in writing down some useful formulas
for the
calculations. Let us consider first the case of the standard one-dimensional
Brownian motion starting from $x\in[-1,1]$ which is killed or
reflected when
hitting the boundaries $-1$ or $1$. We shall write $D$ for Dirichlet condition
on the boundary and $N$ for Neumann condition, which of course
correspond to
killing and, respectively, reflection. Furthermore we shall note, for example,
$p_{\mathit{DN}}(t,x_1,x_2)$ the density of the Brownian motion on $[-1,1]$
killed when
hitting $-1$ and reflected on $1$ more precisely the order in the indices
indicates the boundary condition in $-1$ and $1$, respectively.


\subsubsection{Reflected Brownian motion on $[-1,1]$}

Let $p_{\mathit{NN}}(t,x_1,x_2)$ denote the probability density function of a Brownian
motion at time $t$, starting from $x_1$ and reflected at $-1$ and $1$.
By using the method of images we get the following formula for
the transition density:
\[
p_{\mathit{NN}}(t,x_1,x_2)= \frac{1}{\sqrt{2\pi t}}\sum_{n=-\infty}^{\infty}
\bigl[e ^{-{(x_1-x_2+4n)^2}/({2t})}+e ^ {-{(x_1+x_2+4n+2)^2}/({2t})} \bigr].
\]
The spectral representation of this density writes
\[
p_{\mathit{NN}}(t,x_1,x_2)= {\frac{1}{2}}+{\sum_{n=1}^\infty}
e^{-{n^2\pi^2}/{8}t}
\cos\biggl(\frac{n\pi}{2}(x_1+1) \biggr)\cos\biggl(\frac{n\pi}{2}(x_2+1) \biggr).
\]
These expressions may be found, for example, in \cite{beck92a}.

\subsubsection{Killed Brownian motion on $[-1,1]$}

Let $p_{\mathit{DD}}(t,x_1,x_2)$ denote the probability density function of a Brownian
motion at time $t$, starting from $x_1$ and killed when it exits from
the interval $[-1,1]$. That is,
\[
p_{\mathit{DD}}(t,x_1,x_2)\,dx_2 = \mathbb{P}_{x_1}[B_t\in dx_2; t <\tau_{\mathit{DD}}],
\]
where $\tau_{\mathit{DD}}=\inf\{ t\geq0; B_t\notin[-1,1]\}$.
Then, by the images' method we have
\[
p_{\mathit{DD}}(t,x_1,x_2)={\frac{1}{\sqrt{2\pi t}}\sum_{n=-\infty
}^\infty}
\bigl[e ^ {-{(x_1-x_2+4n)^2}/({2t})}-e ^{-{(x_1+x_2+4n+2)^2}/({2t})}\bigr].
\]
For the law of the exit time we get
\[
\mathbb{P}_{x_1} [\tau_{\mathit{DD}} \in\vd t] ={\frac{1}{\sqrt{2\pi
t^3}}}{\sum_{n=-\infty}^\infty} (-1)^n (x_1+2n+1) e^{-
{(x_1+2n+1)^2}/({2t})}\,\vd t.
\]
The spectral representation can be also written and yields
\[
p_{\mathit{DD}}(t,x_1,x_2)= {\sum_{n=1}^\infty} e ^{-{n^2\pi
^2}/{8}t}\sin\biggl(\frac{n\pi}{2}(x_1+1) \biggr)
\sin\biggl(\frac{n\pi}{2}(x_2+1) \biggr).
\]
The law of the exit time is given by
\[
\mathbb{P}_{x_1} [\tau_{\mathit{DD}} \in\vd t]= \frac{\pi}{2}\sum
_{n=0}^\infty(-1)^n(2n+1)e^ {-{(2n+1)^2\pi^2}/{8}t}\cos\biggl(
\biggl(n+\frac{1}{2} \biggr)\pi x_1 \biggr)\,\vd t.
\]
These expressions may be found, for example, in \cite{beck92a}
or in \cite{milstein99a}.

\subsubsection{Mixed boundary conditions for the Brownian motion on $[-1,1]$}

We give here explicit solutions for the Brownian motion killed on $-1$
and reflected on~$1$.
Let $p_{\mathit{DN}}(t,x_1,x_2)$ denote the probability density function of a Brownian
motion at time $t$, starting from $x_1$ and killed when it hits $-1$
and reflected on~$1$. Then, by the images' method, one gets
\[
p_{\mathit{DN}}(t,x_1,x_2)=\frac{1}{\sqrt{2\pi t}}\sum_{n=-\infty
}^\infty(-1)^n \bigl[e ^{-{(x_1-x_2+4n)^2}/({2t})}-e^{-
{(x_1+x_2+4n+2)^2}/({2t})} \bigr].
\]
Let us denote also by $\tau_{\mathit{DN}}$ the killing time for the Brownian
motion on $[-1,1)$ killed on $-1$ and reflected on $1$. Hence
\[
\mathbb{P}_{x_1} [\tau_{\mathit{DN}} \in\vd t]=\frac{1}{\sqrt{2\pi t^3}}
\sum_{n=-\infty}^\infty(-1)^n (x_1+4n+1)e^{-
{(x_1+4n+1)^2}/({2t})}\,\vd t.
\]
The spectral representation can be also written and yields
\begin{eqnarray*}
p_{\mathit{DN}}(t,x_1,x_2)
&=& \sum_{n=0}^\infty e ^{-{(2n+1)^2\pi^2}/{32}t}
\sin\biggl(\frac{(2n+1)\pi}{4}(x_1+1) \biggr)\\
&&\hspace*{13.7pt}{}\times\sin\biggl(\frac{(2n+1)\pi}{4}(x_2+1) \biggr).
\end{eqnarray*}
Then we get from the spectral representation the law of this exit time,
\[
\mathbb{P}_{x_1} [\tau_{\mathit{DN}} \in\vd t]= \frac{\pi}{8}\sum
_{n=0}^\infty(2n+1)
e^{-{(2n+1)^2\pi^2}/{32}t}\sin\biggl(\frac{(2n+1)\pi}{4}(x_1+1) \biggr)\,\vd t.
\]
The dual situation (reflection on $-1$ and absorption
on $1$) can be obtained easily by the transformation
\[
p_{\mathit{ND}}(t,x_1,x_2)=p_{\mathit{DN}}(t,-x_1,x_2).
\]
These expressions may be found, for example, in \cite{beck92a}.

\subsection{Brownian motion with drift $\mu$}

As in the previous part of the \hyperref[app-1]{Appendix} we consider
here the case of
the Brownian motion with drift on
the interval $[-1,1]$ which is killed or reflected on $-1$ and $1$.
If we note by $p^{L,\mu}_{\cdot\cdot}(t,x_1,x_2)$ the law of the process with
drift $\mu$ and living on $[-L,L]$ and $p^{\mu}_{\cdot\cdot}(t,x_1,x_2)$ the
corresponding law on $[-1,1]$, then by the properties of the Brownian
motion we have
\[
p^{L,\mu}_{\cdot\cdot}(t,x_1,x_2)=\frac{1}{L}p^{\mu L}_{\cdot\cdot} \biggl(\frac
{t}{L^2},\frac{x_1}{L},\frac{x_2}{L} \biggr),
\]
where the dots in the indices can take the value $D$ for a Dirichlet condition
or $N$ for a Neumann condition, as previously noted.

\subsubsection{Brownian motion with drift $\mu$ reflected on $[-1,1]$}

We keep the same notation as before.
The use of the images' method gives the following representation of the density:
\begin{eqnarray*}
p^\mu_{\mathit{NN}}(t,x_1,x_2) &=& \frac{2\mu e^{2\mu x_2}}{e^{2\mu
}-e^{-2\mu}}
+\frac{1}{\sqrt{2\pi t}}\sum_{n=-\infty}^\infty e ^{4\mu
n}e^{-{(x_1-x_2+\mu t+4n)^2}/({2t})}\\
&&{} + \frac{1}{\sqrt{2\pi t}}\sum_{n=-\infty}^\infty
e^{-2\mu x_1}e^{-\mu(4n+2)}e^{-{(x_1+x_2-\mu t +4n +2)^2}/({2t})}\\
&&{} - \mu e^{2\mu x_2}{\sum_{n=-\infty}^\infty} e ^{\mu
(4n+2)}\operatorname{erfc}\biggl(\frac{x_1+x_2+\mu t+4n +2}{\sqrt{2t}} \biggr).
\end{eqnarray*}
This formula can be obtained also from the results in Veestraeten
\cite{veertraeten05a}.

By the spectral method (see, e.g.,
\cite{linetsky05a}), we have, after some calculations,
\begin{eqnarray*}
p^\mu_{\mathit{NN}}(t,x_1,x_2)&=&\frac{2\mu e ^{2\mu x_2}}{e^{2\mu
}-e^{-2\mu}}\\
&&{} + e ^{\mu(x_2-x_1)-{\mu^2}/{2}t}\\
&&\hspace*{9.36pt}{}\times\sum
_{n=1}^\infty\frac{e ^{-{n^2\pi^2}/{8}t}}{\mu^2+
{n^2\pi^2}/{4}}
\biggl[\frac{\pi n}{2}\cos\biggl(\frac{\pi n}{2}(x_1+1) \biggr)\\
&&\hspace*{104pt}{} + \mu
\sin\biggl(\frac{\pi n}{2}(x_1+1) \biggr) \biggr]\\
&&\hspace*{36.45pt}{} \times\biggl[\frac{\pi n}{2}\cos\biggl(\frac{\pi n}{2}(x_2+1) \biggr)+\mu
\sin\biggl(\frac{\pi n}{2}(x_2+1) \biggr) \biggr].
\end{eqnarray*}
%

\subsubsection{Brownian motion with drift $\mu$ on $[-1,1]$ killed at
the boundary}

We keep the same notation as before. By using classical properties of the
Brownian motion and the results from Milstein and Tretyakov \cite{milstein99a}
we have the following transformation:
\[
p^\mu_{\mathit{DD}}(t,x_1,x_2)=e^{\mu(x_2-x_1)-{\mu^2 t}/{2}}p_{\mathit{DD}}(t,x_1,x_2).
\]
Then, by the images' method,
\begin{eqnarray*}
p_{\mathit{DD}}^\mu(t,x_1,x_2) &=& e^{\mu(x_2-x_1)-{\mu^2 t}/{2}}
\frac{1}{\sqrt{2\pi t}}\\
&&{}\times
\sum_{n=-\infty}^\infty\bigl[e^{-
{(x_1-x_2+4n)^2}/({2t})}-e^{-{(x_1+x_2+4n+2)^2}/({2t})} \bigr].
\end{eqnarray*}
We write down both distribution and density for the exit time. The
distribution writes
\begin{eqnarray*}
\mathbb{P}_{x_1} [\tau_{\mathit{DD}}^\mu< t] &=& 1-{\frac{1}{2}} {\sum
_{n=-\infty}^\infty} e ^{4\mu n} \biggl[\operatorname{erfc}\biggl(\frac{x_1+\mu
t +4n
-1}{\sqrt{2t}} \biggr)\\
&&\hspace*{75.3pt}{} -
\operatorname{erfc}\biggl(\frac{x_1+\mu t +4n + 1}{\sqrt{2t}} \biggr) \biggr]\\
&&\hspace*{0pt}{} + {\frac{1}{2}}{\sum_{n=-\infty}^\infty} e ^{-(2\mu x_1 +
\mu(4 n+2))} \biggl[\operatorname{erfc}\biggl(\frac{x_1-\mu t +4n +1}{\sqrt
{2t}} \biggr)\\
&&\hspace*{124.6pt}{} - \operatorname{erfc}\biggl(\frac{x_1-\mu t +4n +
3}{\sqrt{2t}} \biggr) \biggr],
\end{eqnarray*}
while for the density we obtain
\begin{eqnarray*}
\mathbb{P}_{x_1} [\tau_{\mathit{DD}}^\mu\in\vd t] &=&
\frac{e^{-\mu x_1-{\mu^2 t}/{2}}}{\sqrt{2\pi t^3}}\sum
_{n=-\infty}^\infty
\bigl[e^{-\mu} (x_1 +4n+1)e^{ -{(x_1+4n+1)^2}/({2t})} \\
&&\hspace*{84.72pt}{} -e^\mu(x_1 +4n-1)e^{-{(x_1+4n-1)^2}/({2t})}\bigr].
\end{eqnarray*}
The spectral representation can be also written and yields
\begin{eqnarray*}
p_{\mathit{DD}}^\mu(t,x_1,x_2)&=& e^{\mu(x_2-x_1)-{\mu^2 t}/{2}}\\
&&{}\times{\sum
_{n=1}^\infty} e^{-{n^2\pi^2}/{8}t}\sin\biggl(\frac{n\pi}{2}(x_1+1)
\biggr)\sin\biggl(\frac{n\pi}{2}(x_2+1) \biggr).
\end{eqnarray*}
The distribution of the exit time is given by
\begin{eqnarray*}
&&\mathbb{P}_{x_1} [\tau_{\mathit{DD}}^\mu< t] \\
&&\qquad= 1- e ^ {-\mu x_1-{\mu^2t}/{2}} {\sum_{n=1}^\infty}
\bigl(e^{-\mu}- (-1)^ne^\mu\bigr)\frac{2n\pi}{4\mu^2+n^2\pi^2}e ^{-
{n^2\pi^2}/{8}t}\\
&&\hspace*{117.47pt}{}\times\sin\biggl(\frac{n\pi}{2}(x_1+1) \biggr) \\
&&\qquad= 1- e ^ {-\mu x_1-{\mu^2t}/{2}}(e^{-\mu}-e ^\mu) {\sum
_{n=1}^\infty} (-1)^n\frac{n\pi}{\mu^2+n^2\pi^2}e ^{-{n^2\pi
^2}/{2}t}\sin( n\pi x_1)\\
&&\qquad\quad{} - e ^ {-\mu x_1-{\mu^2t}/{2}}(e^{-\mu}+e ^\mu) {\sum
_{n=0}^\infty} (-1)^n\frac{2(2n+1)\pi}{4\mu^2+(2n+1)^2\pi^2}\\
&&\hspace*{162.4pt}{}\times
e^{-{(2n+1)^2\pi^2}/{8}t}\cos\biggl(\frac{(2n+1)\pi}{2}x_1
\biggr)
\end{eqnarray*}
and
\begin{eqnarray*}
\mathbb{P}_{x_1} [\tau_{\mathit{DD}}^\mu\in\vd t]&=&e ^ {-\mu x_1-{\mu^2
t}/{2}}\\
&&\hspace*{0pt}{}\times{ \sum_{n=1}^\infty} \frac{n\pi}{4}\bigl(e^{-\mu}-
(-1)^ne^{\mu}\bigr) e ^{-{n^2\pi^2}/{8}t}\\
&&\hspace*{25.4pt}{}\times\sin\biggl(\frac{n\pi
}{2}(x_1+1) \biggr)\,\vd t.
\end{eqnarray*}
In a more detailed expression we can write this on the form
\begin{eqnarray*}
\mathbb{P}_{x_1} [\tau_{\mathit{DD}}^\mu\in\vd t]&=& e ^ {-\mu x_1-{\mu^2
t}/{2}} (e^{-\mu}-e^{\mu})\\
&&{}\times
{\sum_{n=1}^\infty}(-1)^n \frac{n\pi}{2}
e ^{-{n^2\pi^2}/{2}t}\sin(n\pi x_1)\\
&&{} + e ^ {-\mu x_1-{\mu^2 t}/{2}} (e^{-\mu}+e^{\mu})\\
&&\hspace*{10pt}{}\times{\sum
_{n=0}^\infty} (-1)^n \frac{(2n+1)\pi}{4} e^{-{(2n+1)^2\pi
^2}/{8}t}\\
&&\hspace*{129.8pt}\hspace*{-94.6pt}{}\times\cos\biggl(\frac{(2n+1)\pi}{2}x_1 \biggr) \,\vd t.
\end{eqnarray*}
These expressions may be found, for example, in \cite{beck92a}
or in \cite{milstein99a}.

\subsubsection{Mixed boundary condition for the Brownian motion on
$[-1,1]$ with drift~$\mu$}

The aim is to express some explicit solutions for the Brownian motion
killed on $-1$ and reflected on $1$.
We  solve now the following eigenvalue problem:
\[
\cases{
\frac{1}{2} \varphi'' (x_1) +\mu\varphi' (x_1) =\lambda\varphi(x_1),
\cr
\varphi(-1)=0,\cr
\varphi' (1)=0.}
\]
We can remark first that if $\varphi_\lambda$ is an eigenfunction for
the eigenvalue $\lambda$ for the preceding PDE, then $\lambda$ is negative.

%
%
\begin{table}
\tabcolsep=0pt
\caption{Eigenvalues and eigenfunctions
for the Dirichlet/Neumann problem with~a~constant~transport~term~$\mu
$}\label{table-eigen}
\begin{tabular*}{\tablewidth}{@{\extracolsep{\fill}}lcc@{}}
\hline
$\bolds\mu$ & $\bolds\lambda$ & $\bolds{\varphi_\lambda}$ \\
\hline
$\mu< \frac{1}{2}$ & $\lambda\leq-\frac{\mu^2}{2}$,  & $\frac
{e ^{-\mu x_1}}{\sqrt{2 (1-({\cos^2(2\sqrt{-\mu^2-2\lambda
})})/({2\mu}) )}}\sin(\sqrt{-\mu^2-2\lambda}(x_1+1))$\\
& $\tan(2\sqrt
{-\mu^2-2\lambda}) =\frac{\sqrt{-\mu^2-2\lambda}}{\mu}$\\
[4pt]
$\mu=\frac{1}{2}$ & $-\frac{1}{8}$ & $\frac{\sqrt{3}}{4}e
^{-{x_1}/{2}}(x_1+1)$\\[3pt]
& $\lambda< -\frac{1}{8}$,  &
$\frac{e ^{-{x_1}/{2}}}{\sqrt{2} |\sin(2\sqrt{
{1}/{4}+2\lambda} ) |}\sin(\sqrt{\frac{1}{4}+2\lambda}(1+x_1) )$ \\[4pt]
& $\tan(2\sqrt{ (\frac{1}{4}+2\lambda)}
) = 2\sqrt{ (\frac{1}{4}+2\lambda)}$ & \\[4pt]
$\mu> \frac{1}{2}$ &$\lambda\geq-\frac{\mu^2}{2}$,  & $\frac{e ^{-\mu
x_1}}{\sqrt{{2 \cosh^2 (2\sqrt{\mu^2+2\lambda})}/{\mu
}-1}}\sinh(\sqrt{\mu^2+2\lambda}(x_1+1))$\\[12pt]
& $\tanh(2\sqrt{\mu^2+2\lambda
})=\frac{\sqrt{\mu^2+2\lambda}}{\mu}$ & $\frac{e ^{-\mu
x_1}}{\sqrt{2 (1-{\cos^2(2\sqrt{-\mu^2-2\lambda})}/({2\mu})
)}}\sin(\sqrt{-\mu^2-2\lambda}(x_1+1))$\\
&$\lambda\leq-\frac{\mu^2}{2}$,  & \\
& $\tan(2\sqrt{-\mu^2-2\lambda})
=\frac{\sqrt{-\mu^2-2\lambda}}{\mu}$\\
\hline
\end{tabular*}
\end{table}

We associate with this problem the corresponding second degree equation
and note
$\Delta= \mu^2+2\lambda$.
After a detailed calculus with respect to the sign of $\Delta$ we can express
the countable set of eigenfunctions and eigenvalues with respect to the
possible values
of $\mu$.
There are three different situations, expressed in Table~\ref
{table-eigen} (see, e.g., \cite{pinsky1998}).
The density $p_{\mathit{DN}}(t,x_1,x_2)$ is obtained by
using the spectral expansion
$p_{\mathit{DN}}(t,x_1,x_2)=\sum_{k\geq0}\exp{\lambda_k t}
\varphi_{\lambda_k}(x_1)\varphi_{\lambda_k}(x_2)$,
where $\cdots\leq\lambda_2\leq\lambda_1< \lambda_0$.
The density $q_{\mathit{DN}}(t,x_1)$ of the exit time is also expressed by
\[
\mathbb{P}_{x_1}[\tau_{\mathit{DN}}\in\vd t]/\vd t=
-\sum_{k\geq0}\lambda_k e^{\lambda_k t}
\phi_{\lambda_k}(x_1)\int_{-1}^1 \phi_{\lambda_k}(x_2)\,\vd
x_2.
\]
\end{appendix}

\section*{Acknowledgment}
The authors are grateful to the referee for his helpful remarks and suggestions.


%
\printaddresses

\end{document}